\documentclass[a4paper,11pt]{article}
\usepackage{xcolor}   
\usepackage{verbatim}
\usepackage{amsmath}  
\usepackage{amssymb}   
\usepackage{amsfonts} 
\usepackage{mathrsfs} 
\usepackage{amsthm}  
\usepackage{indentfirst}
\usepackage[title]{appendix}
\usepackage{hyperref} 
\usepackage{graphicx}
\usepackage{setspace}
\hypersetup{hypertex=true,
	colorlinks=true,
	linkcolor=blue,
	anchorcolor=blue,
	citecolor=blue}
\usepackage{microtype} 
\theoremstyle{plain}
\newtheorem{theorem}{Theorem}[section]
\newtheorem{corollary}[theorem]{Corollary}
\newtheorem{proposition}[theorem]{Proposition}
\newtheorem{lemma}[theorem]{Lemma}
\theoremstyle{remark}
\newtheorem{remark}{Remark}
\theoremstyle{definition}
\newtheorem{definition}{Definition}[section]

\numberwithin{equation}{section}
\allowdisplaybreaks[4] 

\begin{document}
\pagenumbering{arabic}

\bigskip\bigskip
\noindent{\Large\bf The Genealogy of  Conditioned Critical Galton--Watson Trees
		\footnote{ This work was supported in part by NSFC (NO. 11971062) and  the National Key Research and Development Program of China (No. 2020YFA0712900).} }
	
\noindent
{Jiayan Guo\footnote{ Email: guojiayan@mail.bnu.edu.cn. School of Mathematical Sciences \& Laboratory of Mathematics and Complex Systems, Beijing Normal University, Beijing 100875, P.R. China.}
		\quad 	
Wenming Hong\footnote{ Email: wmhong@bnu.edu.cn. School of Mathematical Sciences \& Laboratory of Mathematics and Complex Systems, Beijing Normal	University, Beijing 100875, P.R. China.} }

\begin{center}
\begin{minipage}{12cm}
\begin{center}\textbf{Abstract}\end{center}
\footnotesize
Let $T$ be the random family tree associated with the critical Galton--Watson process $\{Z_{n}\}_{n\geq0}$. Geiger (1999) provided an explicit representation of the law of $T$ conditioned on $\{Z_{n}>0\}$ by inductively constructing the conditioned tree $\tilde{T}_n$ along the line of descent of the left--most particle. Intuitively, $\tilde{T}_n$ is composed of a ``spine'', with a conditioned GW tree attached to each of its ``left'' siblings and an ordinary GW tree attached to each of its ``right'' siblings.

The structure of the ``right'' subtree is partially discussed in Geiger (1999), whereas the information about the ``left''  subtree remains largely unexplored. In this paper we investigate the genealogy of the ``left'' and ``right''  subtrees of $\tilde{T}_n$, and some intrinsic properties will be specified,  including the branching mechanism, the limiting behavior of the coalescent times, the limiting number of siblings of the spine, and the scaling limit distribution of the particles. We prove that the sequence of coalescent times of the ``right'' subtree converge in distribution to a sequence of ``nested uniform random variables'' by scaling, interestingly, so do the sequence of coalescent times for the ``left''  subtree after a functional transformation. We also establish the asymptotic independence of the two sequences of coalescent times, which in turn implies the asymptotic independence of the ``left'' and ``right''  subtrees. Finally, as an application, we present a probabilistic proof of the conditional limit theorem for the critical GW process established by Spitzer (unpublished) and Lamperti and Ney (1968), that is, for any fixed $0<t<1$, $Z_{[nt]}/n$ converges in distribution to the sum of two independent exponential random variables with different parameters, which exactly come from the ``left'' and ``right'' subtrees respectively.

\bigskip

\textbf{Keywords:} critical branching process, conditioned tree, reduced process, coalescent time, conditional limit theorem.

\textbf{Mathematics Subject Classification}: 60J80; 60F05
\end{minipage}
\end{center}

\section{Introduction}
Let $\{Z_{n}\}_n$ be a branching process defined by $Z_{0}=1$ and  
	\begin{equation*}
		Z_{n}=\sum_{i=1}^{Z_{n-1}}\xi_{n,i}, \quad n\in\mathbb{N_{+}}, 
	\end{equation*}
	 where $\{\xi_{n,i}\}_{n,i\in\mathbb{N_{+}}}$ is an i.i.d. sequence with distribution $(p_{k})_{k\geq0}$.  Use $\xi$ for the generic copy, $\alpha:=\mathbf{E}\xi$ for the mean, and $f(s):=\mathbf{E}s^{\xi}$ for the generating function, then the $n$th iteration $f_{n}(s)$ is the generating function of $Z_{n}$. When $\alpha<1(=1, >1)$, we say the process is subcritical (critical, supercritical). In this paper we mainly focus on the processes satisfying $\alpha=1$ and $\sigma^{2}=\sum_{k}k(k-1)p_{k}\in(0,\infty)$. 
	 
It is known that the process gets eventually extinct if and only if $\alpha\leq 1$ and $p_{1}\neq1$, so in the critical case, the asymptotic behavior of  $\{Z_{n}\}$ is typically studied under the condition   $\{Z_{n}>0\}$, by analyzing the generating function $f_{n}(s)$ (for example, see Athreya and Ney \cite{Athreya}, Kesten, Ney and Spitzer \cite{Kesten}, Yaglom \cite{Yaglom}). This approach is usually concise and effective. 

Tree structure is also a powerful tool for studying the critical branching processes. Let $T$ be the random family tree associated with $\{Z_{n}\}_{n\geq0}$. We think of $T$ as a rooted planar tree with the distinguishable offspring of each vertex ordered from left to the right.	By inductively constructing the conditioned tree $\tilde{T}_n$ along the line of descent of the left--most particle, Geiger \cite{Geiger99} provided an explicit representation of the law of $T$ conditioned on $\{Z_{n}>0\}$ (Lemma \ref{lem:Geiger} and Theorem \ref {th:Geiger}).

Intuitively, $\tilde{T}_n$ is composed of a ``spine'', with a conditioned GW tree attached to each of its ``left'' siblings and an ordinary GW tree attached to each of its ``right'' siblings. Based on this construction, the representation of distribution of $Z_{m}$ $(0<m\leq n)$ conditioned on $\{Z_{n}>0\}$ can be given (Theorem \ref{th:decom-1}). Using a ``reduced'' method and removing the siblings of the ``spine'' that have no descendant at generation $m$ in $\tilde{T}_n$, we obtain another representation of the conditional law of $Z_{m}$ (Theorem \ref{th:decom-2}). 

Thus, in order to analyze the limiting behavior of the process $\{Z_n\}_{n}$ conditioned on $\{Z_n>0\}$, it is necessary to specify the genealogy and intrinsic properties of Geiger's conditioned tree. Our analysis comprises three  parts: the genealogy of the ``left''   subtree, the genealogy of the ``right''   subtree, and the asymptotic independence of the two subtrees.

The structure of the ``right'' subtree is partially discussed in Geiger (1999), whereas the information about the ``left''  subtree remains largely unexplored. Therefore, we begin by analyzing its branching structure, which is figured out by a time-inhomogeneous branching process (Theorem \ref{th:bpve}). We also prove that: the coalescent times of the left subtree converge in distribution to a sequence of ``nested uniform random variables'' after a functional transformation; for any $0<t<1$, asymptotically there is only one left-sibling of the spine that has descendants at generation $[nt]$; the population size of the left subtree at generation $[nt]$ converges to the exponential distribution after scaling (Theorem \ref{th:MRCA-l}--\ref{th:yaglom-l}). And analogous results are established for the right subtree with $0<t\leq 1$(Theorem \ref{th:MRCA-r}--\ref{th:yaglom-r}).

For the independence, we establish that the two sequences of coalescent times are asymptotically independent (Theorem \ref{th:MRCA-in}). Then as corollaries, we obtain the distribution of generation of the most recent common ancestor of generation $[nt]$ (Corollary \ref{cor:MRCA}), and that the left and right parts at generation $[nt]$ in Geiger's conditioned tree are asymptotically independent (Theorem \ref{th:Z-in}).

As an application, we also present a probabilistic proof of the conditional limit theorem for critical Galton--Watson process, which is established by Spitzer (unpublished, see \cite{Athreya}) and Lamperti and Ney \cite{Lamperti} that for any fixed $0<t<1$, $Z_{[nt]}/n$ converges in distribution to the sum of two independent exponential random variables with different parameters (Theorem \ref{th:Spitzer}). It reveals that the two exponential  random variables are resulted from the ``left'' and ``right'' parts  of the ``spine'' of the Geiger's tree at generation $[nt]$.
 
 \subsection{Geiger's Conditioned Tree} 
 Let $T^{(i)}, 1\leq i\leq Z_{1}(T) $, be the subtrees founded by the children of the root $\emptyset$. Let 
	  \[R_{n+1}(T):=\min\{1\leq i\leq Z_{1}(T): Z_{n}(T^{(i)})>0\}\]
	  denote the rank of the left-most child of the root who has a descendant at generation $n+1$, where $\min\emptyset=\infty$. Then 
	  \begin{lemma}[\cite{Geiger99}]\label{lem:Geiger}
	  	The subtrees $T^{(i)}$, $1\leq i\leq Z_{1}(T)$, are conditionally independent given $\{R_{n+1}(T)=j, Z_{1}(T)=k\}, 1\leq j\leq k<\infty$, with
	  	\begin{equation*}
	  	\mathscr{L}(T^{(i)}|R_{n+1}(T)=j,Z_{1}(T)=k)=\left\{
	  		\begin{aligned}
	  		\begin{array}{ll}
	  		\mathscr{L}(T|Z_{n}=0), &1\leq i\leq j-1;\\
	  		\mathscr{L}(T|Z_{n}>0), &i=j;\\
	  		\mathscr{L}(T), &j+1\leq i\leq k.
	  		\end{array}
	  		\end{aligned}
	  		\right.
	  	\end{equation*}
	  	The conditioned joint distribution of $R_{n+1}(T)$ and $Z_{1}(T)$ $(n\geq0)$ is 
	  	\begin{equation}\label{eq:RZ}
	  		\mathbf{P}(R_{n+1}(T)=j,Z_{1}(T)=k|Z_{n+1}(T)>0)=c_{n}p_{k}\mathbf{P}(Z_{n}=0)^{j-1},
	  	\end{equation}
	  	where $c_{n}=\mathbf{P}(Z_{n}>0)\mathbf{P}(Z_{n+1}>0)^{-1}$.
	  \end{lemma}

\begin{theorem}[Construction of the Conditioned Tree (\cite{Geiger99})]\label {th:Geiger}
	Let $\{(V_{n+1}, W_{n+1} )\}_{n\geq0}$ be a sequence of independent random variables with the distribution \eqref{eq:RZ},
	  	\begin{equation}\label{eq:VW}
	  \mathbf{P}(V_{n+1}=j, W_{n+1}=k)=c_{n}p_{k}\mathbf{P}(Z_{n}=0)^{j-1},\quad1\leq j\leq k<\infty.
	 \end{equation}
 Denote $X_{n+1}:= W_{n+1}- V_{n+1}$. Let $\tilde{T}_{0}$ be an ordinary Galton-Watson tree independent of $(V_{n+1}, W_{n+1})$. Inductively construct $\tilde{T}_{n+1}(n\geq0)$, by the following procedure:
 
	  (i) Let the first generation size of $\tilde{T}_{n+1}$ be $W_{n+1}$.
	  
	  (ii) Take $\tilde{T}_{n}$ to be the subtree founded by the $V_{n+1}$th particle (denoted as $s_{n+1}$) of the first generation of $\tilde{T}_{n+1}$.
	  
	  (iii) Let the $V_{n+1}-1$ siblings to the left of the distinguished first generation particle $s_{n+1}$ found independent Galton-Watson trees conditioned on extinction at generation $n$.
	  
	  (iv) Let the $X_{n+1}$ siblings to the right of the distinguished first generation particle  $s_{n+1}$ found independent ordinary Galton-Watson trees.
	  
	  Then for all $n\geq0$, 
	  \[\mathscr{L}(\tilde{T}_{n})=\mathscr{L}(T|Z_{n}>0).\]\end{theorem}

	\begin{figure}[htbp] 
	\centering
	\includegraphics[scale=0.6]{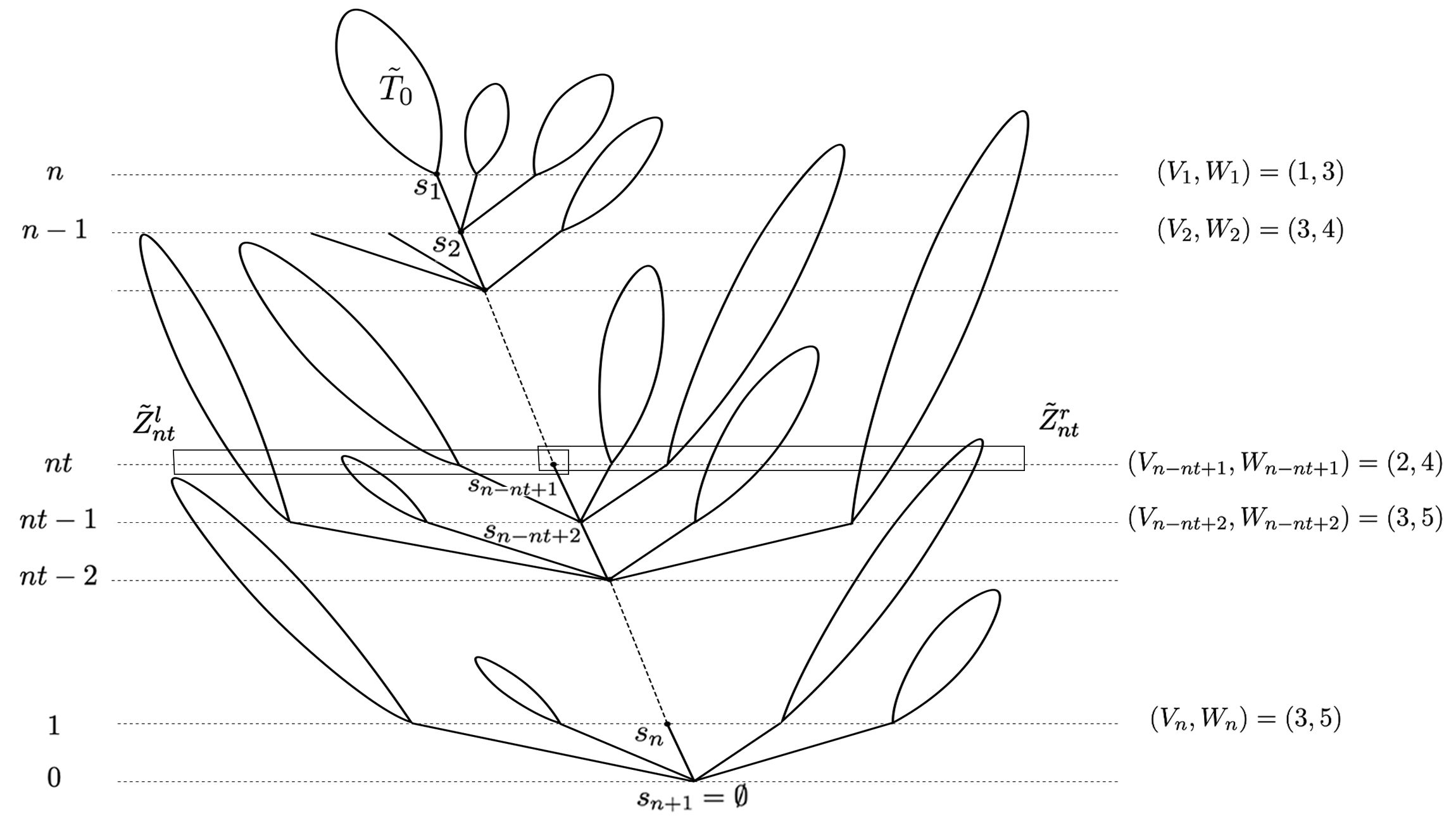} 
	\caption{ {\small Geiger's conditioned tree $\tilde{T}_{n}$. Notice that $\{s_{i}\}$ and $\{(V_{i}, X_{i})\}$ are consistent with the notations in Theorem \ref {th:Geiger}, being in reverse order from $n$ to $1$, while the order of numbering in generation is from $1$ to $n$.}}
	\label{figure1}
	\end{figure}
	
	In other words, $\tilde{T}_{n}$ is a realization for the tree $T$ under the condition $Z_{n}>0$ as follows (see Figure \ref{figure1}): (i) Establish the ``spine'' containing the founding ancestor $\emptyset$ and $n$ distinguished particles $s_{n},\cdots,s_{1}$ at generation $1$ to $n$. (ii) Generate independent random variables $\{(V_{i}, X_{i} )\}_{1\leq i\leq n}$ with notation $X_{i}:= W_{i}-V_{i}$ and distribution \eqref{eq:VW}. Supply $V_{i}-1$ siblings to the left of $s_i$ and $X_{i}$ siblings to the right of $s_i$. (iii) Attach independent trees conditioned to have heights less than $i-1$ to the $V_{i}-1$ siblings to the left of $s_{i}$, and independent unconditioned trees to the $X_{i}$  siblings to the right of $s_{i}$. (iv) Complete the tree by adding an independent, unconditioned tree on the top of the distinguished particle $s_{1}$ at generation $n$.
	
 \subsubsection{Representation of the Conditional Law of $Z_{m}$$(0<m\leq n)$}

	Based on Theorem \ref {th:Geiger}, Geiger \cite{Geiger99} gave the representation of the law of $Z_{n}$ conditioned on $Z_{n}>0$, i.e.,
	\begin{equation}\label{eq:Zn}
		\mathscr{L}(Z_{n}|Z_{n}>0)\overset{\text{d}}{=}1+\sum_{i=1}^{n}\sum_{j=1}^{X_{i}}Z_{i-1;j},
	\end{equation}
	where $Z_{i;j}$ and $X_{i+1}$ ($i\geq0$, $j\geq1$) are independent random variables with $\mathscr{L}(Z_{i;j})=\mathscr{L}(Z_{i})$, and used \eqref{eq:Zn} to  give probabilistic proofs of non-extinction probability and distribution of the most recent common ancestor.
	
	In fact, not only the representation of the conditioned law of $Z_{n}$  can be given, but also that of $Z_{m}(0<m<n)$ . It is worth noting that when considering $Z_{n}$, there is no particle to the left of the spine at generation $n$, so only the properties of the right subtree need to be analyzed. But here when we give the representation of $Z_{m}(0<m<n)$, both particles to the left and right will make contributions. 
	
	Let $\tilde{Z}_m$ be the population size at generation $m$ of the conditioned tree $\tilde{T}_{n}$. Then by Theorem \ref {th:Geiger}, $\mathscr{L}(\tilde{Z}_m)=\mathscr{L}(Z_m|Z_n>0)$, and we can decompose $\tilde{Z}_{m}$ into two parts: the ``left'' part $\tilde{Z}_{m}^{l}$ (population size to the left of the ``spine'') and the ``right'' part $\tilde{Z}_{m}^{r}$ (population size to the right of the ``spine''), and both of them contain the ``spine''. 	
	
	For the right part $\tilde{Z}_{m}^{r}$, the structure is similar to \eqref{eq:Zn}. Standing at generation $m$, let $\tilde{Z}_{m}^{r,0}=1$, i.e., the distinguished particle $s_{n-m+1}$. And for $1\leq i\leq m$, let the increment $\tilde{Z}_{m}^{r,i}-\tilde{Z}_{m}^{r,i-1}$ be the number of particles at generation $m$ produced by the $X_{n-m+i}:=W_{n-m+i}-V_{n-m+i}$ right siblings of $s_{n-m+i}$, then
	\begin{equation}\label{eq:decom-1r}
	 \tilde{Z}_{m}^{r,i}=\tilde{Z}_{m}^{r,i-1}+\sum_{j=1}^{X_{n-m+i}}Z_{m;j}^{r,i-1},
	 \end{equation}
	 where $\{Z_{m;j}^{r,i}\}_{j\geq1}$ are i.i.d. copies of $Z_{m}^{r,i}$ with distribution 
	\begin{equation}\label{eq:decom-1rr}
	\mathscr{L}(Z_{m}^{r,i})=\mathscr{L}({Z}_{i}),
	\end{equation}
	and $\{Z_{n}\}$ is the ordinary Galton--Watson process with $Z_0=1$.

	 For the left part $\tilde{Z}_{m}^{l}$, let $\tilde{Z}_{m}^{l,0}:=1$, i.e., the distinguished particle $s_{n-m+1}$. And for $1\leq i\leq m$, let the increment $\tilde{Z}_{m}^{l,i}-\tilde{Z}_{m}^{l,i-1}$ be the number of particles at generation $m$ produced by the $(V_{n-m+i}-1)$ left siblings of $s_{n-m+i}$, then
	\begin{equation}\label{eq:decom-1l}
	 \tilde{Z}_{m}^{l,i}=\tilde{Z}_{m}^{l,i-1}+\sum_{j=1}^{V_{n-m+i}-1}Z_{m;j}^{l,i-1},
	 \end{equation}
	 where  $\{Z_{m;j}^{l,i}\}_{j\geq1}$ are i.i.d. copies of $Z_{m}^{l,i}$ with distribution 
	\begin{equation}\label{eq:decom-1ll}
	\mathscr{L}(Z_{m}^{l,i})=\mathscr{L}({Z}_{i}|Z_{i+n-m}=0).
	\end{equation}

	Summarize the above with Theorem \ref {th:Geiger} and denote $\tilde{Z}_{m}^{l}:=\tilde{Z}_{m}^{l,m}$, $\tilde{Z}_{m}^{r}:=\tilde{Z}_{m}^{r,m}$, we get the following theorem, which is consistent with \eqref{eq:Zn} when $m=n$.
	\begin{theorem}\label{th:decom-1}
	For all $n\geq0$, $0<m\leq n$,
	\[\mathscr{L}(\tilde{Z}_{m}^{l},\tilde{Z}_{m}^{r})=\mathscr{L}\left(1+\sum_{i=1}^{m}\sum_{j=1}^{V_{n-m+i}-1}Z_{m;j}^{l,i-1},1+\sum_{i=1}^{m}\sum_{j=1}^{X_{n-m+i}}Z_{m;j}^{r,i-1}\right),\]
	where $\{Z_{m;j}^{l,i}\}_{i,j}$ and $\{Z_{m;j}^{r,i}\}_{i,j}$ are independent. In particular,
	\[\mathscr{L}(Z_m|Z_n>0)=\mathscr{L}(\tilde{Z}_m)=\mathscr{L}(\tilde{Z}_m^l+\tilde{Z}_m^r-1)\overset{\text{d}}{=}1+\sum_{i=1}^{m}\left(\sum_{j=1}^{V_{n-m+i}-1}Z_{m;j}^{l,i-1}+\sum_{j=1}^{X_{n-m+i}}Z_{m;j}^{r,i-1}\right).\]
	\end{theorem}
	
\subsubsection{Reduced Geiger's Conditioned Tree}
\begin{figure}[htbp]
	\centering
	\includegraphics[scale=0.7]{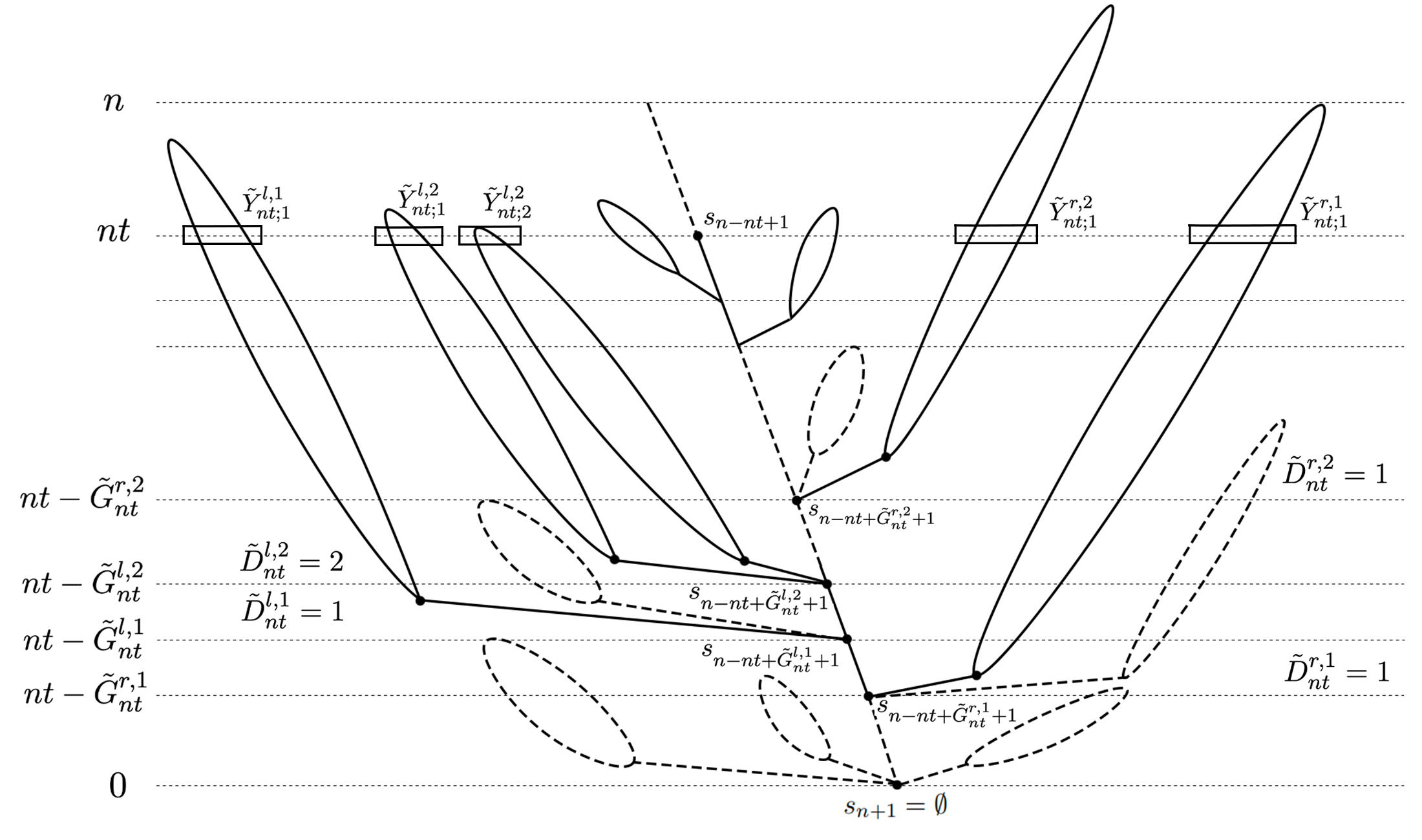}
	\caption{ {\small The reduced Geiger's tree for generation $m$. Only particles having non-empty descendants at $m$ are left. }}
	\label{figure2}
\end{figure}

	Since we are interested in the particles at generation $m$,  only the particles having non-empty descendants at $m$ in Geiger's tree $\tilde{T}_n$ will make contributions. Therefore, using the ``reduced'' method, we can remove the siblings of $s_{i}$  ($n-m< i\leq n$, recall that here is in reverse order, see Figure \ref{figure2}) that have no descendant at $m$ in $\tilde{T}_n$, and decompose the particles at generation $m$ in another way. 

To do this, we first introduce the sequence of coalescent times $\{\tilde{G}_{m}^{l,k}\}_{k}$ for the left subtree. Intuitively, by \eqref{eq:decom-1l}, $\tilde{Z}_{m}^{l}$ is composed of a spine particle and $m$ blocks, where the $i$th block is produced by the $(V_{n-m+i}-1)$ left siblings of $s_{n-m+i}$, and we are concerned with the distributions of the most common ancestors of these $(m+1)$ components. 

For $0<m\leq n$, denote
\[\tilde{G}_{m}^{l}=\tilde{G}_{m}^{l,1}:=\min\left\{0\leq i\leq m |  \tilde{Z}_{m}^{l,i}=\tilde{Z}_{m}^{l,m}\right\},\]
then  $(m-\tilde{G}_{m}^{l})$ is the generation of the most recent common ancestor of  all the $(m+1)$ blocks on the left subtree. Furthermore, for $k\geq2$, define
\[\tilde{G}_{m}^{l,k}:=\min\left\{0\leq i\leq \tilde{G}_{m}^{l,k-1}-1 | \tilde{Z}_{m}^{l,i}=\tilde{Z}_{m}^{l,\tilde{G}_{m}^{l,k-1}-1}\right\},\]
where $\tilde{G}_{m}^{l,k}=0$ if $\tilde{G}_{m}^{l,k-1}=0$. Then $(m-\tilde{G}_{m}^{l,k})$ is the generation of the most recent common ancestor of the remaining components on the left subtree after removing the first $(k-1)$ non-empty blocks from the outside inward, i.e., the most recent common ancestor of particles that are on the left part of generation $m$ of $\tilde{T}_n$  and rooted by $s_{n-m+\tilde{G}_{m}^{l,k-1}}$.

Then we have 
\[m\geq \tilde{G}_{m}^{l,1}\geq\tilde{G}_{m}^{l,2}\geq\cdots \geq\tilde{G}_{m}^{l,m}\geq0=\tilde{G}_{m}^{l,m+1}=\cdots,\]
where $\tilde{G}_{m}^{l,k}=\tilde{G}_{m}^{l,k+1}$ holds if and only if $\tilde{G}_{m}^{l,k}=\tilde{G}_{m}^{l,k+1}=0$. 

Define the sequence of coalescent times $\{\tilde{G}_{m}^{r,k}\}_{k}$ for the right subtree in the same way.

For $k\geq1$, denote $\tilde{H}_{m}^{l,k}$ for the number of particles  satisfying: (i) is on the left part of generation $m$ of $\tilde{T}_n$; (ii) rooted by $s_{n-m+\tilde{G}_{m}^{l,k}}$. Then by \eqref{eq:decom-1l} and the construction of Geiger's tree, for $0\leq p\leq m-1$,
\begin{equation}\label{eq-tl}
	\mathscr{L}\left(\tilde{H}_{m}^{l,k}\big|\tilde{G}_{m}^{l,k}=p+1\right)
	=\mathscr{L}\left(1+\sum_{i=1}^{p}\sum_{j=1}^{V_{n-m+i}-1}Z_{m;j}^{l,i-1}\right).
\end{equation}

Denote $\tilde{D}_{m}^{l,k}$ for the number of particles satisfying: (i) is a sibling to the left  of $s_{n-m+\tilde{G}_{m}^{l,k}}$; (ii) has non-empty descendants at generation $m$ of $\tilde{T}_n$. And use $\tilde{Y}_{m;j}^{l,k}$ to represent the number of descendants at generation $m$ produced by the $j$th of the $\tilde{D}_{m}^{l,k}$ particles. Write $\tilde{Y}_{m}^{l,k}$ for the generic copy, then by \eqref{eq:decom-1ll}, for $0\leq p\leq m-1$,
 \begin{equation}\label{eq-yl}
 	\mathscr{L}\left(\tilde{Y}_{m}^{l,k}\big|\tilde{G}_{m}^{l,k}=p+1\right)
 	=\mathscr{L}\left(Z_{m}^{l,p}\big|Z_{m}^{l,p}>0\right)
 	=\mathscr{L}\left(Z_{p}\big|Z_{p}>0,Z_{p+n-m}=0\right),
 \end{equation}
 where to avoid separating the case we define $\tilde{H}_{m}^{l,k}=1$, $\tilde{D}_{m}^{l,k}=\tilde{Y}_{m;1}^{l,k}=0$ if $\tilde{G}_{m}^{l,k}=0$.

Use the similar notations $\tilde{H}_{m}^{r,k}$, $\tilde{D}_{m}^{r,k}$, $\tilde{Y}_{m;j}^{r,k}$ for the right subtree. Then we have, for $0\leq p\leq m-1$,
\begin{equation}\label{eq-tr}
	\mathscr{L}\left(\tilde{H}_{m}^{r,k}\big|\tilde{G}_{m}^{r,k}=p+1\right)
	=\mathscr{L}\left(1+\sum_{i=1}^{p}\sum_{j=1}^{X_{n-m+i}}Z_{m;j}^{r,i-1}\right),
\end{equation}
and
 \begin{equation}\label{eq-yr}
 	\mathscr{L}\left(\tilde{Y}_{m}^{r,k}\big|\tilde{G}_{m}^{r,k}=p+1\right)
 	=\mathscr{L}\left(Z_{m}^{r,p}\big|Z_{m}^{r,p}>0\right)
 	=\mathscr{L}\left(Z_{p}\big|Z_{p}>0\right).
 \end{equation}
Similarly, define $\tilde{H}_{m}^{r,k}=1$, $\tilde{D}_{m}^{r,k}=\tilde{Y}_{m;1}^{r,k}=0$ if $\tilde{G}_{m}^{r,k}=0$.

The detailed distributions of \eqref{eq-tl}-\eqref{eq-yr} will be discussed in Section \ref{section:2} and Section \ref{section:3}. We also note that by the independence of subtrees and branching property, given $\tilde{G}_{m}^{l,k}$, $\{\tilde{Y}_{m;j}^{l,k}\}_j$ are i.i.d. and independent of $\tilde{H}_{m}^{l,k}$. Given $\tilde{G}_{m}^{r,k}$, $\{\tilde{Y}_{m;j}^{r,k}\}_j$ are i.i.d. and independent of $\tilde{H}_{m}^{r,k}$. Given $\tilde{G}_{m}^{r,k}$ and $\tilde{G}_{m}^{r,k}$, $\{\tilde{Y}_{m;j}^{l,k}\}_j$ and $\{\tilde{Y}_{m;j}^{r,k}\}_j$ are independent.

In summary, we get a reduced version representation for the conditional law of $Z_{m}$.
\begin{theorem}\label{th:decom-2}
	For all $n\geq0$, $0<m\leq n$, $k_l,k_r\geq 1$,
	\[\mathscr{L}(\tilde{Z}_{m}^{l},\tilde{Z}_{m}^{r})=\mathscr{L}\left(\tilde{H}_{m}^{l,k_l}+\sum_{i=1}^{k_l}\sum_{j=1}^{\tilde{D}_{m}^{l,i}}\tilde{Y}_{m;j}^{l,i},\tilde{H}_{m}^{r,k_r}+\sum_{i=1}^{k_r}\sum_{j=1}^{\tilde{D}_{m}^{r,i}}\tilde{Y}_{m;j}^{r,i}\right),\]
	where $\{\tilde{Y}_{m;j}^{l,k}\}_{k,j}$ and $\{\tilde{Y}_{m;j}^{r,k}\}_{k,j}$ are independent given $\{\tilde{G}_{m}^{l,k}\}_{k}$ and $\{\tilde{G}_{m}^{r,k}\}_{k}$. In particular,
	\[\mathscr{L}(\tilde{Z}_{m})\overset{\text{d}}{=}1+\sum_{i=1}^{m}\left(\sum_{j=1}^{\tilde{D}_{m}^{l,i}}\tilde{Y}_{m;j}^{l,i}+\sum_{j=1}^{\tilde{D}_{m}^{r,i}}\tilde{Y}_{m;j}^{r,i}\right).\]
\end{theorem}

\subsection{Main Results}
In this subsection, we investigate the genealogy of the left and right subtrees, encompassing the coalescent times $\{\tilde{G}_{m}^{l,k}\}_{k}$, $\{\tilde{G}_{m}^{r,k}\}_{k}$, the number of siblings $\{\tilde{D}_{m}^{l,k}\}_{k}$, $\{\tilde{D}_{m}^{r,k}\}_{k}$, the limiting distributions and asymptotic independence of $\tilde{Z}_{m}^{l}$ and $\tilde{Z}_{m}^{r}$.

\subsubsection{Genealogy of the Right Subtree}
For the coalescent times of the right subtree at generation $[nt]$, we prove that they converge in distribution to the sequence of ``nested uniform random variables''. For simplicity we drop the square bracket, but $nt$ will always be understood as $[nt]$.
\begin{definition}[nested uniform random variables]\label{def:nurv}
If $U_{1}$ is a uniform random variable on $[a,b]$, and inductively, $U_{k}$ $(k\geq2)$ is uniform on $[a,U_{k-1}]$, then $\{U_{k}\}_{k\geq1}$ is said to be a sequence of {\it nested uniform random variables} on $[a,b]$.
\end{definition}
\begin{theorem}\label{th:MRCA-r}
Suppose $\alpha=1$ and $\sigma^{2}<\infty$, then for any fixed  $0<t\leq 1$ and $k\geq 1$, as $n\rightarrow\infty$, \[\frac{\tilde{G}_{nt}^{r,k}}{nt}\overset{\text{d}}{\rightarrow}U_{k},\] where $\{U_{k}\}_{k\geq1}$ is a sequence of nested uniform random variables on $[0,1]$.
\end{theorem}

For the number of siblings, we prove that asymptotically there is only one right-sibling that has descendants at generation $[nt]$.  
\begin{theorem}\label{th:D-r}
	Suppose $\alpha=1$ and $\sigma^{2}<\infty$, then for any fixed $0<t\leq1$ and $k\geq1$, as $n\rightarrow\infty$,
	\[\tilde{D}_{nt}^{r,k}\overset{\text{d}}{\rightarrow}1.\]
\end{theorem}

We also obtain the Yaglom's Theorem for the right subtree.
\begin{theorem}\label{th:yaglom-r}
	Suppose $\alpha=1$ and $\sigma^{2}<\infty$, then for any fixed $0<t\leq1$, as $n\rightarrow\infty$,
	\[\mathbf{P}\left(\frac{\tilde{Z}_{nt}^{r}}{n}\geq x\right)\rightarrow \exp\left(-\frac{2x}{t\sigma^2}\right), \quad x\geq 0.\]
\end{theorem}

\begin{remark} \label{remark:1}
	As mentioned in \eqref{eq:Zn}, when $t=1$, $\mathscr{L}(\tilde{Z}^{r}_{n})=\mathscr{L}(Z_{n}|Z_n>0)$, and the above theorems are consistent with the known results about the critical Galton--Watson process. For example, Theorem \ref{th:MRCA-r} $(k=1)$ is that conditioned on $Z_{n}>0$, the generation of the most recent common ancestor of generation $n$ is asymptotically uniformly distributed, which is due to Zubkov \cite{Zubkov} and proved by Geiger \cite{Geiger99} in a probabilistic way; Theorem \ref{th:D-r} $(k=1)$ is contained in Geiger \cite{Geiger00} that given non-extinction at generation $n$ there are asymptotically no more than two children of the most recent common ancestor with a descendant at $n$; Theorem \ref{th:yaglom-r} is the Yaglom's Theorem (Kesten et al. \cite{Kesten}, Yaglom \cite{Yaglom}) that conditioned on $Z_{n}>0$, $Z_{n}/n$ converges to the exponential distribution with parameter $2/\sigma^{2}$. Here we extend the results to $0<t<1$ and $k>1$.
\end{remark}

\subsubsection{Genealogy of the Left Subtree}

For the limiting behavior of the left subtree, we begin by analyzing its branching mechanism. Notice that by \eqref{eq:decom-1rr}, the branching mechanism of the right subtree is the normal. However, for the left part, it is a new situation and is heavily influenced  by the condition $\{Z_{n}=0\}$. As a result, \eqref{eq:decom-1ll} is figured out by a time-inhomogeneous branching process (branching process in varying environment).
\begin{theorem}\label{th:bpve}
	For fixed $n\geq0$, the process $\{\hat{Z}_{k}^{(n)}\}_{k\geq0}$ defined by
	\[\mathscr{L}(\hat{Z}_{k}^{(n)}):=\mathscr{L}({Z}_{k}|Z_{n}=0)\]
	is a  branching process in varying environment $\hat{Z}_{k}^{(n)}=\sum_{j=1}^{\hat{Z}_{k-1}^{(n)}}\hat{\xi}^{(n)}_{k,j},$
	where the generating function of offspring is 
	\begin{align*}
	\mathbf{E}s^{\hat{\xi}_{k}^{(n)}}=\left\{
	\begin{aligned}
		&\frac{f(f_{n-k}(0)s)}{f_{n-k+1}(0)},\quad & k<n; \\
		&1,\quad & k\geq n.
	\end{aligned}
	\right.
\end{align*}
\end{theorem}

And for the coalescent times of the left subtree at generation $[nt]$, interestingly, we prove that they still converge in distribution to a sequence of ``nested uniform random variables'' after a functional transformation, although the behavior of the left subtree is heavily influenced by the condition $\{Z_{n}=0\}$. 
\begin{theorem}\label{th:MRCA-l}
Suppose $\alpha=1$ and $\sigma^{2}<\infty$, then for any fixed $0<t<1$ and $k\geq1$, as $n\rightarrow\infty$, \[g_{t}\left(\frac{\tilde{G}_{nt}^{l,k}}{nt}\right)\overset{\text{d}}{\rightarrow}U_{k},\] where $g_{t}(x):=x(tx+1-t)^{-1}$ and $\{U_{k}\}_{k\geq1}$ is a sequence of nested uniform random variables on $[0,1]$.
\end{theorem}

Notice that since there is no particle to the left of the spine at generation $n$ in the conditioned tree $\tilde{T}_{n}$, here the $t$ must be strictly less than $1$. 

For the number of siblings, we also have,

\begin{theorem}\label{th:D-l}
	Suppose $\alpha=1$ and $\sigma^{2}<\infty$, then for any fixed $0<t<1$ and $ k\geq 1$, as $n$ $\rightarrow\infty$,
	\[\tilde{D}_{nt}^{l,k}\overset{\text{d}}{\rightarrow}1.\]
\end{theorem}

We also obtain the Yaglom's Theorem  of the left subtree.
\begin{theorem}\label{th:yaglom-l}
	Suppose $\alpha=1$ and $\sigma^{2}<\infty$, then for any fixed $0<t<1$, as $n\rightarrow\infty$,
	\[\mathbf{P}\left(\frac{\tilde{Z}_{nt}^{l}}{n}\geq x\right)\rightarrow \exp\left(-\frac{2x}{t(1-t)\sigma^2}\right), \quad x\geq 0.\]
\end{theorem}

\begin{remark}
	 Fleischmann and Siegmund-Schultze \cite{Fleischmann76} considered the reduced Galton--Watson process $(\mu^{t}(t'))_{0\leq t'\leq t}$, where $\mu^{t}(t')$ is the number of individuals which are alive at generation $t'$ and have descendants alive at generation $t$ in the family tree of an ordinary Galton-Watson process $(\mu(t))_{t}$. Under the assumptions of criticality and finite variance, they proved that the reduced process converges to a ``transformed Yule process''. That is, a continuous time Markov branching process on $[0,\infty)$ with initial distribution $\delta_{1}$, where each particle arising at time $\varepsilon$ has a lifetime uniformly distributed on $[\varepsilon,1)$ and at each split exactly two particles arise. 	
	This property of ``two branches'' and 	``uniform lifetime'' also appeared in invariance principle for reduced family trees of critical spatially homogeneous branching processes and conditional central limit theorem for critical branching random walk, see Fleischmann and Siegmund-Schultze \cite{Fleischmann77} and Hong and Liang \cite{Hong}. 
	
	Here, Geiger's conditioned tree method analyzes critical Galton-Watson processes from a different perspective: a ``spine'' divides the tree into  a ``left" subtree and a ``right" subtree. Our results reveal that the ``right'' subtree continues to exhibit the``two branches'' and ``uniform lifetime'' properties (Theorem \ref{th:MRCA-r}--\ref{th:D-r}), while the ``left'' subtree is heavily influenced by the condition $\{Z_{n}=0\}$ (Theorem\ref{th:MRCA-l}).
\end{remark}

\subsubsection{Asymptotic Independence of the Two Subtrees}
From  Theorem \ref{th:decom-1}, it is evident that $\tilde{Z}_{nt}^{l}$ and $\tilde{Z}_{nt}^{r}$ are heavily dependent since $V_{n-nt+i}$ and $X_{n-nt+i}$ satisfies the joint distribution in \eqref{eq:VW} for each $i$, although $\{Z_{nt;j}^{l,i}\}_{i,j}$ and $\{Z_{nt;j}^{r,i}\}_{i,j}$ are independent. However, we prove that the two sequences of coalescent times are asymptotically independent.
\begin{theorem}\label{th:MRCA-in}
	Suppose $\alpha=1$ and $\sigma^{2}<\infty$, then for any fixed $0<t<1$, $k_{l},k_{r}\geq1$, and $x,y\in[0,t]$,
	\[\lim_{n\rightarrow\infty}\mathbf{P}(\tilde{G}_{nt}^{l,k_{l}}\leq nx, \tilde{G}_{nt}^{r,k_{r}}\leq ny)
	=\lim_{n\rightarrow\infty}\mathbf{P}(\tilde{G}_{nt}^{l,k_{l}}\leq nx)
	\cdot\lim_{n\rightarrow\infty}\mathbf{P}(\tilde{G}_{nt}^{r,k_{r}}\leq ny).\]
\end{theorem}

Thus as a corollary, let $G_{nt}$ be the generation of the most recent common ancestor of all particles at generation $[nt]$ of tree $T$, we get
\begin{corollary}\label{cor:MRCA}
	Suppose $\alpha=1$ and $\sigma^{2}<\infty$, then for any fixed $0<t<1$,
	\[\lim_{n\rightarrow\infty}\mathbf{P}\left(\frac{G_{nt}}{n}\leq x\Big|Z_{n}>0\right)=1-\frac{(t-x)^2}{t^2(1-x)}, \quad 0\leq x\leq t.\]
\end{corollary} 
For more results about the coalescence time of some or all particles at generation $n$, we refer to \cite{Athreya12, Athreya16}, \cite{HarrisS}, \cite{Zubkov}.

And after using the ``reduced'' method, it directly follows from Theorem \ref{th:decom-2}, \ref{th:D-r},  \ref{th:D-l} and  \ref{th:MRCA-in} that the left and right parts at generation $[nt]$ in Geiger's conditioned tree are asymptotically independent.
\begin{theorem}\label{th:Z-in}
	Suppose $\alpha=1$ and $\sigma^{2}<\infty$, then for any fixed $0<t<1$ and $x,y\in[0,t]$,
	\[\lim_{n\rightarrow\infty}\mathbf{P}(\tilde{Z}_{nt}^{l}\leq nx, \tilde{Z}_{nt}^{r}\leq ny)
	=\lim_{n\rightarrow\infty}\mathbf{P}(\tilde{Z}_{nt}^{l}\leq nx)
	\cdot\lim_{n\rightarrow\infty}\mathbf{P}(\tilde{Z}_{nt}^{r}\leq ny).\]
\end{theorem}

\subsection{Application: Probabilistic Proof of a Conditional Limit Theorem}
  Spitzer (unpublished, see \cite{Athreya}), Lamperti and Ney \cite{Lamperti} proved that 
  \begin{theorem}\label{th:Spitzer}
	Suppose $\alpha=1$ and $\sigma^{2}<\infty$, then for any fixed $0<t<1$, as $n\rightarrow\infty$,
	\[\mathscr{L}\left(\frac{Z_{nt}}{n}\Big|Z_{n}>0\right)\overset{\text{d}}{\rightarrow}U_{t}+V_{t},\]	where $U_{t}$ and $V_{t}$ are independent random variables having exponential distributions with parameters $2/(t(1-t)\sigma^{2})$ and $2/(t\sigma^{2})$ respectively. 
\end{theorem}

The proof in \cite{Athreya} is short and elegant based on the Laplace transform. However, this method cannot visually explain the origin of these two exponential distributions. Here based on the genealogy of Geiger's conditioned tree, we give a probabilistic proof for Theorem \ref{th:Spitzer} and reveals that the two  exponential random variables $U_{t}$ and $V_{t}$ come from the ``left'' and ``right'' subtrees of the ``spine'' of Geiger's tree at generation $[nt]$ conditioned on non-extinction at generation $n$, respectively. We also note that Zhang et al. \cite{Zhang} have done some related work on Theorem \ref{th:Spitzer}. They roughly provided a characterization of $V_t$, but that of $U_t$ and the asymptotic independence were not considered.

	As mentioned in Remark \ref{remark:1}, when $t=1$, only right particles remain at generation $n$ in the Geiger's conditioned tree, and thus the result reduces to the classical Yaglom's limit. For the Yaglom's Theorem, Lyons, Pemantle and Peres \cite{Lyons} gave a probabilistic proof using the size-biased Galton--Watson tree; Geiger \cite{Geiger99,Geiger00} proved the Yaglom's Theorem by  Geiger's conditioned tree and the distribution of most recent common ancestor of particles at generation $n$; Ren, Song and Sun \cite{Ren,Ren2} proposed a two-spine decomposition of the critical Galton--Watson tree and used it to give another probabilistic proof of the theorem; Cardona-Tob\'on and Palau \cite{Cardona} gave the probabilistic proof when environment is varying. 
	
	Here we consider the case where $t$ is strictly less than one. In this regime, both the left and right subtrees contribute, and thus the genealogy and intrinsic properties of Geiger's conditioned tree must be investigated.
\

The article is organized as follows. In Section \ref{section:2} we analyze the genealogy of the left part, which includes the branching mechanism and moment, the limiting behavior for the coalescent times, the limiting number of siblings of the spine, and the limiting distribution, and present the proofs of Theorem \ref{th:bpve}--\ref{th:yaglom-l}. Section \ref{section:3} addresses parallel questions for the right part, providing the proofs of Theorem \ref{th:MRCA-r}--\ref{th:yaglom-r}. In Section \ref{section:4} we establish the asymptotic independence for the two sequences of coalescent times, which implies the independence of the left and right subtrees; this section contains the proofs  of Theorem \ref{th:MRCA-in}--\ref{th:Z-in}. Finally, Section \ref{section:5} is devoted to the proof of Theorem \ref{th:Spitzer}.

\section{Genealogy of the Left Subtree}\label{section:2}
\subsection{Branching Mechanism and Moments}
\begin{proof}[Proof of Theorem \ref{th:bpve}]
For $j\geq 0$,
\begin{align*}
	\mathbf{P}(\hat{Z}^{(n)}_{k}=j)=\left\{
	\begin{aligned}
	&\frac{\mathbf{P}({Z}_{k}=j)[\mathbf{P}(Z_{n-k}=0)]^{j}}{\mathbf{P}(Z_{n}=0)},&\quad k<n;\\
	&\mathbf{1}_{\{j=0\}},&\quad k\geq n. 
	\end{aligned}
	\right.
\end{align*}

To verify the branching property, recall that $f(s)$ is the generating function of $\xi$, the $n$th iteration $f_{n}(s)$ is the generating function of $Z_{n}$ (with $Z_{0}=1$), and $[f_{n}(s)]^{l}$ is the generating function of $Z_{n}$ (with $Z_{0}=l$).

$(i)$ When $k< n$, for all $j,l\in\mathbb{N}$, 
\begin{align*}	
\mathbf{P}(\hat{Z}^{(n)}_{k}=j|\hat{Z}^{(n)}_{k-1}=l)
&=\frac{\mathbf{P}(Z_{k}=j, Z_{k-1}=l|Z_{n}=0)}{\mathbf{P}(Z_{k-1}=l|Z_{n}=0)}\\
&=\frac{\mathbf{P}(Z_{k}=j, Z_{n}=0|Z_{k-1}=l)}{\mathbf{P}(Z_{n}=0|Z_{k-1}=l)}\\
&=\mathbf{P}(Z_{1}=j|Z_{0}=l)\frac{[\mathbf{P}(Z_{n-k}=0)]^{j}}{[\mathbf{P}(Z_{n-k+1}=0)]^{l}}.
\end{align*}	

Thus for all $0<s<1$,
\begin{align}\label{eq-bpve1}
	\sum_{j=0}^{\infty}\mathbf{P}(\hat{Z}^{(n)}_{k}=j|\hat{Z}^{(n)}_{k-1}=l)s^{j}
	&=\sum_{j=0}^{\infty}\mathbf{P}(Z_{1}=j|Z_{0}=l)\frac{[\mathbf{P}(Z_{n-k}=0)]^{j}}{[\mathbf{P}(Z_{n-k+1}=0)]^{l}}s^{j}\notag\\
	&=\left(\frac{f(f_{n-k}(0)s)}{f_{n-k+1}(0)}\right)^{l}\notag\\
	&=\left(\sum_{j=0}^{\infty}\mathbf{P}(Z_{1}=j|Z_{0}=1)\frac{[\mathbf{P}(Z_{n-k}=0)]^{j}}{\mathbf{P}(Z_{n-k+1}=0)}s^{j}\right)^{l}\notag\\
	&=\left(\sum_{j=0}^{\infty}\mathbf{P}(\hat{Z}^{(n)}_{k}=j|\hat{Z}^{(n)}_{k-1}=1)s^{j}\right)^{l}.
\end{align}

$(ii)$ When $k=n$, for all $j,l\in\mathbb{N}$,
\[\mathbf{P}(\hat{Z}^{(n)}_{k}=j|\hat{Z}^{(n)}_{k-1}=l)=\mathbf{1}_{\{j=0\}},\]
thus
\begin{equation}\label{eq-bpve2}
\sum_{j=0}^{\infty}\mathbf{P}(\hat{Z}^{(n)}_{k}=j|\hat{Z}^{(n)}_{k-1}=l)s^{j}=1=\left(\sum_{j=0}^{\infty}\mathbf{P}(\hat{Z}^{(n)}_{k}=j|\hat{Z}^{(n)}_{k-1}=1)s^{j}\right)^{l}.
\end{equation}

$(iii)$ When $k>n$, 
\begin{equation}\label{eq-bpve3}
	\hat{Z}^{(n)}_{k-1}=\hat{Z}^{(n)}_{k}=0.
\end{equation}

Combining \eqref{eq-bpve1}--\eqref{eq-bpve3}, we get that $\{\hat{Z}_{k}^{(n)}\}_{k\geq0}$ is a branching process in varying environment defined by 
\[\hat{Z}_{k}^{(n)}=\sum_{j=1}^{\hat{Z}_{k-1}^{(n)}}\hat{\xi}^{(n)}_{k,j},\]
where $\{\hat{\xi}^{(n)}_{k,j}\}_{j\geq1}$ is an i.i.d. sequence with the distribution 
\begin{align*}
	\mathbf{P}(\hat{\xi}^{(n)}_{k}=j)=\left\{
	\begin{aligned}
		&\frac{p_{j}\cdot[\mathbf{P}(Z_{n-k}=0)]^{j}}{\mathbf{P}(Z_{n-k+1}=0)},\quad & k<n; \\
		&\mathbf{1}_{\{j=0\}},\quad & k\geq n,
	\end{aligned}
	\right.
\end{align*}
and generating function 
\begin{align*}
	\hat{g}_{k}^{(n)}(s):=\mathbf{E}s^{\hat{\xi}_{k}^{(n)}}=\left\{
	\begin{aligned}
		&\frac{f(f_{n-k}(0)s)}{f_{n-k+1}(0)},\quad & k<n; \\
		&1,\quad & k\geq n.
	\end{aligned}
	\right.
\end{align*}
\end{proof}
Recall that by \eqref{eq:decom-1ll}, $Z_{nt}^{l,i}=\hat{Z}_{i}^{(i+n-nt)}$. In the subsequent proof, the moments of $Z_{nt}^{l,i}$ are needed, and we have the following propositions. 
\begin{proposition}\label{prop:moment-l1}
	Suppose $\alpha=1$ and $\sigma^{2}<\infty$, then as $n\rightarrow\infty$, for any $0\leq i\leq nt$,
	\[\mathbf{E}Z_{nt}^{l,i}\sim\left(\frac{n(1-t)}{i+n(1-t)}\right)^{2}.\]
\end{proposition}
\begin{proof}
	For any $0\leq k\leq i+n-nt$,
	\[\frac{\mathrm{d}\hat{g}_{k}^{(i+n-nt)}(s)}{\mathrm{d}s}\Big|_{s=1}
	=\frac{f'(f_{i+n-nt-k}(0))f_{i+n-nt-k}(0)}{f_{i+n-nt-k+1}(0)},\]
	and thus
	\[\mathbf{E}(Z_{nt}^{l,i})
	=\prod_{k=1}^{i}\frac{\mathrm{d}\hat{g}_{k}^{(i+n-nt)}(s)}{\mathrm{d}s}\Big|_{s=1}
	=\frac{f_{n-nt}(0)}{f_{i+n-nt}(0)}\cdot\frac{f'_{i+n-nt}(0)}{f'_{n-nt}(0)}.\]
	
	Since $1-f_{n}(0)=\mathbf{P}(Z_{n}>0)\sim 2/(\sigma^{2}n)$ and there exist positive constant $c$, s.t., $\lim_{n\rightarrow\infty}n^{2}f'_{n}(0)=c$ (see corollary 9.4 of \cite{Athreya} for example), then the theorem follows. 
\end{proof}
\begin{proposition}\label{prop:moment-l2}
	Suppose $\alpha=1$ and $\sigma^{2}<\infty$, then as $n\rightarrow\infty$, for any $0\leq i\leq nt$,
	\[\mathbf{E}[Z_{nt}^{l,i}(Z_{nt}^{l,i}-1)]\sim i\left(\frac{n(1-t)}{i+n(1-t)}\right)^{3}\sigma^{2}.\]
\end{proposition}
\begin{proof}
	For $0\leq k\leq  i+n-nt$,
	\[\frac{\mathrm{d}^{2}\hat{g}_{k}^{(i+n-nt)}(s)}{\mathrm{d}s^{2}}\Big|_{s=1}=\frac{f''(f_{i+n-nt-k}(0))[f_{i+n-nt-k}(0)]^{2}}{f_{i+n-nt-k+1}(0)},\]
	then by Lemma 4 of Kersting \cite{Kersting} we have, for the branching process in varying environment,  
	\begin{align*}
	\frac{\mathbf{E}[Z_{nt}^{l,i}(Z_{nt}^{l,i}-1)]}{(\mathbf{E}Z_{nt}^{l,i})^{2}}
	&=\sum_{k=1}^{i}\frac{\mathrm{d}^{2}\hat{g}_{k}^{(i+n-nt)}(s)}{\mathrm{d}s^{2}}\Big|_{s=1}
	\cdot\left(\frac{\mathrm{d}\hat{g}_{k}^{(i+n-nt)}(s)}{\mathrm{d}s}\Big|_{s=1}\right)^{-2}
	\cdot\frac{1}{\mathbf{E}\hat{Z}_{k-1}^{(i+n-nt)}},
	\end{align*}
	so by $f'(1)=1$ and $f''(1)=\sigma^2$, as $n\rightarrow\infty$,
	\begin{align*}
		\mathbf{E}[Z_{nt}^{l,i}(Z_{nt}^{l,i}-1)]
		&=\frac{f^2_{n-nt}(0)}{f_{i+n-nt}(0)}\sum_{k=1}^{i}\frac{f''(f_{i+n-nt-k}(0))f'_{i+n-nt-k+1}(0)f'_{i+n-nt}(0)}{[f'(f_{i+n-nt-k}(0))]^2[f'_{n-nt}(0)]^2}\\
		&\sim\sigma^{2}\frac{(n-nt)^{4}}{(i+n-nt)^{2}}\sum_{k=1}^{i}\frac{1}{(i+n-nt-k+1)^{2}}\\
		&\sim i\left(\frac{n-nt}{i+n-nt}\right)^{3}\sigma^{2},
	\end{align*}
	where the last line is by the fact that 
	\[\lim_{n\rightarrow\infty}n\sum_{k=n}^{\infty}k^{-2}=\lim_{n\rightarrow\infty}n^{-1}\sum_{k=n}^{\infty}(k/n)^{-2}=\int_{1}^{\infty}x^{-2}\mathrm{d}x=1.\]
\end{proof}

By Theorem \ref{th:decom-1},
\begin{equation*}
	\tilde{Z}_{nt}^{l}\overset{\text{d}}{=}1+\sum_{j=1}^{V_{n-nt+1}-1}Z_{nt;j}^{l,0}+\sum_{j=1}^{V_{n-nt+2}-1}Z_{nt;j}^{l,1}+\cdots+\sum_{j=1}^{V_{n}-1}Z_{nt;j}^{l,nt-1},
\end{equation*}
thus by Proposition \ref{prop:moment-l1} and \ref{prop:moment-l2} we can estimate the moments of $\tilde{Z}_{nt}^{l}$.
\begin{theorem}\label{th:moment-l1}
	Suppose $\alpha=1$ and $\sigma^{2}<\infty$, then for any fixed $0<t<1$,
	\[\lim_{n\rightarrow\infty}n^{-1}\mathbf{E}{\tilde{Z}_{nt}^{l}}=\frac{t(1-t)\sigma^{2}}{2}.\]
\end{theorem}
\begin{proof}
	We just need to show that
	\[\lim_{n\rightarrow\infty}\frac{\mathbf{E}{(\tilde{Z}_{nt}^{l}-1)}}{nt}=\lim_{n\rightarrow\infty}\frac{\sum_{i=1}^{nt}\mathbf{E}(V_{n-nt+i}-1)\mathbf{E}Z_{nt}^{l,i-1}}{nt}=\frac{(1-t)\sigma^{2}}{2}.\]
	
	Since by the dominated convergence theorem $\mathbf{E}(V_{n}-1)\rightarrow \mathbf{E}(V_{\infty}-1)=\sigma^{2}/2$, we have
	\[\lim_{n\rightarrow\infty}\frac{\sum_{i=1}^{nt}|\mathbf{E}(V_{n-nt+i}-1)-2^{-1}\sigma^{2}|}{nt}=0,\]
	so if
	\begin{equation}\label{eq2-1}
		\lim_{n\rightarrow\infty}\max_{1\leq i\leq nt}\mathbf{E}Z_{nt}^{l,i-1}<\infty
	\end{equation}
	and
	\begin{equation}\label{eq2-2}
		\lim_{n\rightarrow\infty}\frac{\sum_{i=1}^{nt}\mathbf{E}Z_{nt}^{l,i-1}}{nt}=1-t
	\end{equation}
	are true, then
	\begin{equation*}
		\begin{aligned}
			&\lim_{n\rightarrow\infty}\left|\frac{\sum_{i=1}^{nt}\mathbf{E}(V_{n-nt+i}-1)\mathbf{E}Z_{nt}^{l,i-1}}{nt}-\frac{\sigma^{2}}{2}\frac{\sum_{i=1}^{nt}\mathbf{E}Z_{nt}^{l,i-1}}{nt}\right|\\
			\leq&\lim_{n\rightarrow\infty}\frac{\sum_{i=1}^{nt}|\mathbf{E}(V_{n-nt+i}-1)-2^{-1}\sigma^{2}|\cdot \mathbf{E}Z_{nt}^{l,i-1}}{nt}\\
			\leq&\lim_{n\rightarrow\infty}\max_{1\leq i\leq nt}\mathbf{E}Z_{nt}^{l,i-1}\cdot\lim_{n\rightarrow\infty}\frac{\sum_{i=1}^{nt}|\mathbf{E}(V_{n-nt+i}-1)-2^{-1}\sigma^{2}|}{nt}=0,
		\end{aligned}
	\end{equation*}
	and the theorem follows.
	
	For \eqref{eq2-1}, notice that for $0\leq i\leq nt-1$,
	\[
			\mathbf{E}Z_{nt}^{l,i}=\frac{f_{n-nt}(0)f'_{i+n-nt}(0)}{f'_{n-nt}(0)f_{i+n-nt}(0)}
	\]
	is decreasing with $i$, so 
	\[\lim_{n\rightarrow\infty}\max_{1\leq i\leq nt}\mathbf{E}Z_{nt}^{l,i-1}
	=\mathbf{E}Z_{nt}^{l,0}=1.\]
	
	For \eqref{eq2-2}, by Proposition \ref{prop:moment-l1},
	\begin{equation*}
	\lim_{n\rightarrow\infty}\frac{\sum_{i=1}^{nt}\mathbf{E}Z_{nt}^{l,i-1}}{nt}
	=\lim_{n\rightarrow\infty}\frac{1}{nt}\sum_{i=1}^{nt}\left(\frac{n(1-t)}{i+n(1-t)}\right)^{2}
	=\int_{0}^{1}\left(1+\frac{t}{1-t}x\right)^{-2}dx
	=1-t.	
	\end{equation*}
\end{proof}
	\begin{theorem}\label{th:moment-l2}
	Suppose $\alpha=1$ and $\sigma^{2}<\infty$, then for any fixed $0<t<1$,
	\[\lim_{n\rightarrow\infty}n^{-2}\mathbf{E}(\tilde{Z}_{nt}^{l})^{2}=\frac{[t(1-t)\sigma^{2}]^{2}}{2}.\]
\end{theorem}
\begin{proof}
By Theorem \ref{th:decom-1},
	\begin{align*}
		&\mathbf{E}(\tilde{Z}_{nt}^{l}-1)^{2}\\	=&\sum_{1\leq k\neq l\leq nt}\mathbf{E}\left(\sum_{j=1}^{V_{n-nt+k}-1}Z_{nt;j}^{l,k-1}\cdot\sum_{j=1}^{V_{n-nt+l}-1}Z_{nt;j}^{l,l-1}\right)
		+\sum_{i=1}^{nt}\mathbf{E}\left(\sum_{j=1}^{V_{n-nt+i}-1}Z_{nt;j}^{l,i-1}\right)^{2}\\
		=&\sum_{1\leq k\neq l\leq nt}\left(\mathbf{E}V_{n-nt+k}-1\right)\cdot \left(\mathbf{E}V_{n-nt+l}-1\right)\cdot \mathbf{E}Z_{nt}^{l,k-1}\cdot \mathbf{E}Z_{nt}^{l,l-1}\\
		&+\sum_{i=1}^{nt}\mathbf{E}(V_{n-nt+i}-1)^{2}\cdot\left(\mathbf{E}Z_{nt}^{l,i-1}\right)^{2}\\
		&+\sum_{i=1}^{nt}\left(\mathbf{E}V_{n-nt+i}-1\right )\cdot \left[\mathbf{E}(Z_{nt}^{l,i-1})^{2}-\left(\mathbf{E}Z_{nt}^{l,i-1}\right)^{2}\right]\\
		=&:I_{n,1}+I_{n,2}+I_{n,3},
	\end{align*}
and we will show that the main order comes from $I_{n,1}$ and $I_{n,3}$.

Using the same method as in the last theorem, by Proposition \ref{prop:moment-l1} and \ref{prop:moment-l2},
\begin{align*}
	\lim_{n\rightarrow\infty}\frac{I_{n,1}}{n^{2}}
	=&\lim_{n\rightarrow\infty}\frac{1}{n^2}\sum_{1\leq k,l\leq nt}\left(\mathbf{E}V_{n-nt+k}-1\right)\cdot \left(\mathbf{E}V_{n-nt+l}-1\right)\cdot \mathbf{E}Z_{nt}^{l,k-1}\cdot \mathbf{E}Z_{nt}^{l,l-1}\\
	=&\left(\frac{\sigma^2}{2}t\right)^2\int_{0}^{1}\int_{0}^{1}\left(1+\frac{t}{1-t}x\right)^{-2}\left(1+\frac{t}{1-t}y\right)^{-2}\mathrm{d}x\mathrm{d}y\\
	=&\left(\frac{\sigma^2}{2}t(1-t)\right)^2,\\
	\lim_{n\rightarrow\infty}\frac{I_{n,3}}{n^{2}}
	=&\lim_{n\rightarrow\infty}\frac{1}{n^{2}}\sum_{i=1}^{nt}\left(\mathbf{E}V_{n-nt+i}-1\right )\cdot \mathbf{E}\left[Z_{nt}^{l,i-1}\left(Z_{nt}^{l,i-1}-1\right)\right]\\
	=&\frac{\sigma^2}{2}t^2\cdot\lim_{n\rightarrow\infty}\frac{1}{nt}\sum_{i=1}^{nt}\frac{i-1}{nt}\left(1+\frac{t}{1-t}\frac{i-1}{nt}\right)^{-3}\sigma^{2}\\
	=&\frac{(\sigma^2t)^2}{2}\int_{0}^{1}x\left(1+\frac{t}{1-t}x\right)^{-3}\mathrm{d}x\\
	=&\left(\frac{\sigma^2}{2}t(1-t)\right)^2.
\end{align*}

As for the rate of $I_{n,2}$, recall that $\mathbf{P}(V_{n}-1=j)=\sum_{i=j+1}^{\infty}c_{n-1}p_{i}\mathbf{P}(Z_{n-1}=0)^{j}$, then by $\sigma^2<\infty$, $\sum_{i=n}^{\infty}i^2p_{i}\rightarrow0$ and $n[1-\mathbf{P}(Z_{n-1}=0)]\rightarrow2/\sigma^2$, we have 
\begin{align*}
	\lim_{n\rightarrow\infty}\frac{\mathbf{E}(V_{n}-1)^2}{n}
	&=\lim_{n\rightarrow\infty}\sum_{j=1}^{\infty}\sum_{i=j+1}^{\infty}j^2p_{i}\frac{\mathbf{P}(Z_{n-1}=0)^{j}}{n}\\
	&\leq\lim_{n\rightarrow\infty}\sum_{j=1}^{\infty}\frac{\mathbf{P}(Z_{n-1}=0)^{j}\sum_{i=j+1}^{\infty}i^2p_{i}}{n}\rightarrow0,
\end{align*}
thus using Proposition \ref{prop:moment-l1} again,
\begin{equation*}
\begin{aligned}
	\lim_{n\rightarrow\infty}\frac{I_{n,2}}{n^{2}}
	=&\lim_{n\rightarrow\infty}\frac{1}{n^{2}}\sum_{i=1}^{nt}\mathbf{E}(V_{n-nt+i}-1)^{2}\cdot\left(\mathbf{E}Z_{nt}^{l,i-1}\right)^{2}\\
	\leq&\lim_{n\rightarrow\infty}\frac{\max_{1\leq i\leq nt}\mathbf{E}(V_{n-nt+i}-1)^{2}}{n}\cdot\lim_{n\rightarrow\infty}\frac{\sum_{i=1}^{nt}\left(\mathbf{E}Z_{nt}^{l,i-1}\right)^{2}}{n}\\
	=&\lim_{n\rightarrow\infty}\frac{\max_{1\leq i\leq nt}\mathbf{E}(V_{n-nt+i}-1)^{2}}{n}\cdot t\int_{0}^{1}\left(1+\frac{t}{1-t}x\right)^{-4}dx
\end{aligned}
\end{equation*}
equals $0$, and the theorem follows.
\end{proof}

\subsection{Coalescent Times}
Recall the coalescent time sequence of the left subtree that
\[\tilde{G}_{nt}^{l}=\tilde{G}_{nt}^{l,1}:=\min\{0\leq i\leq nt | \tilde{Z}_{nt}^{l,i}=\tilde{Z}_{nt}^{l,nt}\},\]
and for $k\geq2$,
\[\tilde{G}_{nt}^{l,k}:=\min\left\{0\leq i\leq \tilde{G}_{nt}^{l,k-1}-1 | \tilde{Z}_{nt}^{l,i}=\tilde{Z}_{nt}^{l,\tilde{G}_{nt}^{l,k-1}-1}\right\}.\]

So when analyzing the distributions for $\{\tilde{G}_{nt}^{l,k}\}_{k}$ and $\{\tilde{G}_{nt}^{r,k}\}_{k}$, the following lemma is needed.
\begin{lemma}
	Suppose $\alpha=1$ and $\sigma^{2}<\infty$, then as $n\rightarrow\infty$, for all $nx\leq i\leq nt-1$,
	\begin{align}
	&\text{(i)}\,\mathbf{P}\left(\tilde{Z}_{nt}^{l,i+1}>\tilde{Z}_{nt}^{l,i}\right)\sim\frac{1}{i}-\frac{1}{i+n-nt},\label{eq3-1}\\
	&\text{(ii)}\,\mathbf{P}\left(\tilde{Z}_{nt}^{l,i+1}>\tilde{Z}_{nt}^{l,i}, \tilde{Z}_{nt}^{r,i+1}=\tilde{Z}_{nt}^{r,i}\right)\sim\frac{1}{i}-\frac{1}{i+n-nt},\label{eq3-2}\\
	&\text{(iii)}\,\mathbf{P}\left(\tilde{Z}_{nt}^{r,i+1}>\tilde{Z}_{nt}^{r,i}\right)\sim\frac{1}{i},\label{eq3-3}\\
	&\text{(iv)}\,\mathbf{P}\left(\tilde{Z}_{nt}^{l,i+1}=\tilde{Z}_{nt}^{l,i}, \tilde{Z}_{nt}^{r,i+1}>\tilde{Z}_{nt}^{r,i}\right)\sim\frac{1}{i},\label{eq3-4}\\
	&\text{(v)}\,1-\mathbf{P}\left(\tilde{Z}_{nt}^{l,i+1}=\tilde{Z}_{nt}^{l,i}, \tilde{Z}_{nt}^{r,i+1}=\tilde{Z}_{nt}^{r,i}\right)
	\sim\frac{2}{i}-\frac{1}{i+n-nt}.\label{eq3-5}
	\end{align}
\end{lemma} 
\begin{proof}
	$(i)$ For $nx\leq i\leq nt-1$, define
	\[J^{(n)}_{i+1}:=\sum_{j=1}^{V_{n-nt+i+1}-1}\mathbf{1}\{Z_{nt;j}^{l,i}>0\},\]
	then by \eqref{eq:decom-1l},
	\[\mathbf{P}(J^{(n)}_{i+1}=1)
	\leq \mathbf{P}(J^{(n)}_{i+1}>0)
	=\mathbf{P}(\tilde{Z}_{nt}^{l,i+1}>\tilde{Z}_{nt}^{l,i})
	\leq \mathbf{E}(J^{(n)}_{i+1}).\]
	
Since 		
	\begin{align*}
	&\mathbf{P}(J^{(n)}_{i+1}=1)=\sum_{k=2}^{\infty}(k-1)\mathbf{P}(V_{n-nt+i+1}=k)\mathbf{P}(Z_{nt}^{l,i}>0)[\mathbf{P}(Z_{nt}^{l,i}=0)]^{k-2},\\
	&\mathbf{E}(J^{(n)}_{i+1})=\mathbf{E}(V_{n-nt+i+1}-1)\mathbf{P}(Z_{nt}^{l,i}>0),\\
	&\mathbf{P}(Z_{nt}^{l,i}=0)=\frac{\mathbf{P}(Z_{i}=0)}{\mathbf{P}(Z_{i+n-nt}=0)},
	\end{align*}
	where $\mathbf{P}(Z_{n}>0)\sim 2/(\sigma^{2}n)$ and $\mathbf{E}(V_n-1)\rightarrow\sigma^{2}/2$,  we have, as $n\rightarrow\infty$, for all $nx\leq i\leq nt-1$,
	\begin{equation*}
		\mathbf{P}(\tilde{Z}_{nt}^{l,i+1}>\tilde{Z}_{nt}^{l,i})\sim \frac{\sigma^{2}}{2}\mathbf{P}(Z_{nt}^{l,i}>0)\sim\frac{1}{i}-\frac{1}{i+n-nt}.
	\end{equation*}
	
$(ii)$ Using the same method,
	\begin{align*}
		&\mathbf{P}\left(\tilde{Z}_{nt}^{l,i+1}>\tilde{Z}_{nt}^{l,i}, \tilde{Z}_{nt}^{r,i+1}=\tilde{Z}_{nt}^{r,i}\right)\\
		=&\mathbf{P}\left(\sum_{k=1}^{V_{n-nt+i+1}-1}Z_{nt;k}^{l,i}>0, \sum_{j=1}^{X_{n-nt+i+1}}Z_{nt;j}^{r,i}=0\right)\\
		=&\sum_{p=1}^{\infty}\sum_{q=0}^{\infty}
		\mathbf{P}\left(V_{n-nt+i+1}-1=p,X_{n-nt+i+1}=q,\sum_{k=1}^{p}Z_{nt;k}^{l,i}>0, \sum_{j=1}^{q}Z_{nt;j}^{r,i}=0 \right)\\
		=&\sum_{p=1}^{\infty}\sum_{q=0}^{\infty}
		\mathbf{P}\left(V_{n-nt+i+1}-1=p,X_{n-nt+i+1}=q\right)
		\mathbf{P}\left(\sum_{k=1}^{p}Z_{nt;k}^{l,i}>0\right)[\mathbf{P}(Z_{nt}^{r,i}=0)]^{q}\\
		\sim&\sum_{p=1}^{\infty}\sum_{q=0}^{\infty}
		\mathbf{P}\left(V_{n-nt+i+1}-1=p,X_{n-nt+i+1}=q\right)p\mathbf{P}(Z_{nt}^{l,i}>0)\\
		=&\mathbf{E}(V_{n-nt+i+1}-1)\mathbf{P}(Z_{nt}^{l,i}>0),
	\end{align*}
	thus as $n\rightarrow\infty$, for all $nx\leq i\leq nt-1$,
	\begin{equation*}
		\mathbf{P}\left(\tilde{Z}_{nt}^{l,i+1}>\tilde{Z}_{nt}^{l,i}, \tilde{Z}_{nt}^{r,i+1}=\tilde{Z}_{nt}^{r,i}\right)\sim \frac{\sigma^{2}}{2}\mathbf{P}(Z_{nt}^{l,i}>0)\sim\frac{1}{i}-\frac{1}{i+n-nt}.
	\end{equation*}
	
	$(iii) (iv)$ Similarly, for all $nx\leq i\leq nt-1$, consider 
	\[J^{(n)}_{i+1}:=\sum_{j=1}^{X_{n-nt+i+1}}\mathbf{1}\{Z_{nt;j}^{r,i}>0\}\]
	to get
	\[
		\mathbf{P}(\tilde{Z}_{nt}^{r,i+1}>\tilde{Z}_{nt}^{r,i})\sim \mathbf{P}\left(\tilde{Z}_{nt}^{l,i+1}=\tilde{Z}_{nt}^{l,i}, \tilde{Z}_{nt}^{r,i+1}>\tilde{Z}_{nt}^{r,i}\right)\sim\frac{\sigma^{2}}{2}\mathbf{P}(Z_{nt}^{r,i}>0)\sim\frac{1}{i}.
	\]

	$(v)$ Just combine \eqref{eq3-1}--\eqref{eq3-4} or consider
	\[J^{(n)}_{i+1}:=
	\sum_{k=1}^{V_{n-nt+i+1}-1}\mathbf{1}\{Z_{nt;k}^{l,i}>0\}
	+\sum_{j=1}^{X_{n-nt+i+1}}\mathbf{1}\{Z_{nt;j}^{r,i}>0\}.\]
	\end{proof}
	
\begin{proof}[Proof of Theorem \ref{th:MRCA-l}]
Let $\{U_{k}\}_{k\geq1}$ be a nested uniform random variable sequence on $[0,1]$ defined in Definition \ref{def:nurv}. Since $(\tilde{Z}_{nt}^{l,i})_{i\geq0}$ have non-negative independent increments, by \eqref{eq3-1},
\begin{align*}
	\mathbf{P}\left(\tilde{G}_{nt}^{l,1}\leq nx\right)
	&=\prod_{i=nx}^{nt-1}\mathbf{P}\left(\tilde{Z}_{nt}^{l,i}=\tilde{Z}_{nt}^{l,i+1}\right)\\
	&\sim\prod_{i=nx}^{nt-1}\left[1-\left(\frac{1}{i}-\frac{1}{i+n-nt}\right)\right]\\
	&\sim \exp\left\{-\sum_{i=nx}^{nt-1}\left(\frac{1}{i}-\frac{1}{i+n-nt}\right)\right\}\\
	&\sim\exp\left\{-\log\frac{nt-1}{nx-1}+\log\frac{nt-1+n-nt}{nx-1+n-nt}\right\}\\
	&\rightarrow\frac{x}{t(x+1-t)}
\end{align*}
for $0\leq x\leq t$. Thus as $n\rightarrow\infty$, for all $0\leq x \leq 1$,
\[\mathbf{P}\left(g_{t}\left(\frac{\tilde{G}_{nt}^{l,1}}{nt}\right)\leq x\right)=\mathbf{P}\left(\frac{\tilde{G}_{nt}^{l,1}}{n}\leq \frac{t(1-t)x}{1-tx}\right)\rightarrow x,\]
i.e., $g_{t}(\tilde{G}_{nt}^{l,1}/nt)\overset{\text{d}}{\rightarrow}U_{1}$.

Assume that for some $k\geq1$, $g_{t}(\tilde{G}_{nt}^{l,k}/nt)\overset{\text{d}}{\rightarrow}U_{k}$. Then for $0\leq x\leq t$, if $ny\leq nx$,
\[\mathbf{P}\left(\tilde{G}_{nt}^{l,k+1}\leq nx\Big|\tilde{G}_{nt}^{l,k}=ny\right)=1;\]
if $ny>nx$,  
\[\mathbf{P}\left(\tilde{G}_{nt}^{l,k+1}\leq nx\Big|\tilde{G}_{nt}^{l,k}=ny\right)\sim \prod_{i=nx}^{ny-2}\left[1-\left(\frac{1}{i}-\frac{1}{i+n-nt}\right)\right]\rightarrow\frac{x(y+1-t)}{y(x+1-t)},\]
thus
\[\mathbf{P}\left(g_{t}\left(\frac{\tilde{G}_{nt}^{l,k+1}}{nt}\right)\leq x\Big|g_{t}\left(\frac{\tilde{G}_{nt}^{l,k}}{nt}\right)=y\right)
\rightarrow\frac{x}{y}.\]

So as $n\rightarrow\infty$, $g(\tilde{G}_{nt}^{l,k+1}/nt)\overset{\text{d}}{\rightarrow}U[0,U_{k}]=U_{k+1}$, and the theorem follows by induction.
\end{proof}

Thus by combining Theorem \ref{th:MRCA-l} with Proposition \ref{prop:nurv}, we get the limiting distribution functions  of the coalescent time sequences of the left subtree.
\begin{corollary}
	Suppose $\alpha=1$ and $\sigma^{2}<\infty$, then for any fixed $0<t<1$, $k\in\mathbb{N}_{+}$, and $0\leq x\leq t$,
\[\lim_{n\rightarrow\infty}\mathbf{P}\left(\frac{\tilde{G}_{nt}^{l,k}}{n}\leq x\right)=\frac{x}{t(x+1-t)}\sum_{m=0}^{k-1}\frac{1}{m!}\left(\ln\frac{t(x+1-t)}{x}\right)^{m}.\]
\end{corollary}

\subsection{Siblings of the Spine}
By Theorem \ref{th:decom-2}, in order to obtain the limiting distribution of siblings $\{\tilde{D}_{nt}^{l,k}\}_{k}$, the distributions of random variables in \eqref{eq-tl}--\eqref{eq-yl} are needed.

Notice that in \eqref{eq-tl}, if $p=ns$ for some fixed $0<s<t<1$, then $1+\sum_{i=1}^{ns}\sum_{j=1}^{V_{n-nt+i}-1}Z_{nt;j}^{l,i-1}$
has the same structure as $\mathscr{L}(\tilde{Z}_{nt}^{l})$ in Theorem \ref{th:decom-2}, where $n(1-t+s)$ is the ``new $n$'' and $s(1-t+s)^{-1}$ is the ``new $t$''. For this reason, we sometimes use the notation $\tilde{Z}_{n,t}^{l}$  to precisely represent the length $n$ and proportion $t$.

For the limiting distributions of $\mathscr{L}\left(Z_{p}\big|Z_{p}>0,Z_{p+n-nt}=0\right)$ in \eqref{eq-yl}, since $p$ will be usually regarded as $\tilde{G}_{nt}^{l,k}$ which converges after scaling, we define $L_{n,t}(n\geq 0, 0<t<1)$ as a random variable with distribution
\[\mathscr{L}(L_{n,t})=\mathscr{L}(Z_{nt}|Z_{nt}>0,Z_{n}=0).\]
	The limiting distributions and estimation of moments for $L_{n,t}$ are given in Lemma \ref{lem:Lnt}.

The main idea of our proof is to obtain expression for the distribution of the standardized random variable
\begin{equation}\label{eq:dotZ-defl}
	\dot{Z}_{nt}^{l}:=\frac{\tilde{Z}_{nt}^{l}}{a_{nt}}\quad\text{and}\quad\dot{L}_{nt}:=\frac{L_{nt}}{b_{nt}},
\end{equation}
where $a_{nt}:=\mathbf{E}\tilde{Z}_{nt}^{l}$, $b_{nt}:=\mathbf{E}L_{nt}$ are the mean values. Then  $\mathbf{E}{\dot{Z}_{nt}^{l}}=\mathbf{E}{\dot{L}_{nt}}=1$, and by Theorem \ref{th:moment-l1}, Lemma \ref{lem:Lnt}, for any fixed $0<t<1$,
\[\lim_{n\rightarrow\infty}\frac{a_{nt}}{n}=\lim_{n\rightarrow\infty}\frac{b_{nt}}{n}=\frac{t(1-t)\sigma^{2}}{2}.\]

\begin{theorem}\label{th:decom-3l}
	For all $n\geq0$ and $0<t<1$, $k\geq1$,
\begin{equation}\label{eq:dotZ-disl}
	\dot{Z}_{nt}^{l}\overset{\text{d}}{=}\frac{a_{N_{nt}^{l,k}T_{nt}^{l,k}}}{a_{nt}}\dot{Z}_{N_{nt}^{l,k}T_{nt}^{l,k}}^{l}+\sum_{i=1}^{k}\frac{b_{N_{nt}^{l,i}T_{nt}^{l,i}}}{a_{nt}}\sum_{j=1}^{\tilde{D}_{nt}^{l,i}}\dot{L}_{N_{nt}^{l,i}T_{nt}^{l,i};j},
\end{equation}
where $N_{nt}^{l,k}:=\tilde{G}_{nt}^{l,k}-1+n-nt$, $T_{nt}^{l,k}:=(\tilde{G}_{nt}^{l,k}-1)/(\tilde{G}_{nt}^{l,k}-1+n-nt)$. And for all $m\geq0$, $0<s<1$, $\{\dot{L}_{m,s;j}\}_{j\geq1}$ are independent copies of $\dot{L}_{m,s}$; $\{\dot{L}_{m,s;j},\dot{Z}_{m,s}^{l}\}_{m,s,j}$ are independent and independent of $\{\tilde{D}_{nt}^{l,k},\tilde{G}_{nt}^{l,k}\}_{k}$.
\end{theorem}
\begin{proof}
	First we prove the case $k=1$. By Theorem \ref{th:decom-2}, we have, for all $n\geq0$ and $0<t<1$,	
	\[\tilde{Z}_{m}^{l}\overset{\text{d}}{=}\tilde{H}_{nt}^{l,1}+\sum_{j=1}^{\tilde{D}_{nt}^{l,1}}\tilde{Y}_{nt;j}^{l,1},\]
	where by the independence of subtrees, given $\tilde{G}_{nt}^{l,1}$, $\{\tilde{Y}_{nt;j}^{l,1}\}_j$ are i.i.d. and independent of $\tilde{H}_{nt}^{l,1}$.
	
	By \eqref{eq-tl}--\eqref{eq-yl}, for $0\leq m\leq nt-1$,
	\begin{align*}
	&\mathscr{L}\left(\tilde{H}_{nt}^{l,1}\big|\tilde{G}_{nt}^{l,1}=m+1\right)=\mathscr{L}\left(\tilde{Z}_{m+n-nt,m/m+n-nt}^{l}\right),\\
	&\mathscr{L}\left(\tilde{Y}_{nt;j}^{l,1}\big|\tilde{G}_{nt}^{l,1}=m+1\right)=\mathscr{L}\left(L_{m+n-nt,m/m+n-nt}\right).
	\end{align*}
	
	Then for all $x\geq0$, by the independence of $\{\dot{Z}^{l}_{i,s}, \dot{L}_{i,s;j}\}_{i,s,j}$ and $(\tilde{D}_{nt}^{l,1},\tilde{G}_{nt}^{l,1})$,
	\begin{align*}
	&\mathbf{P}\left(\dot{Z}_{nt}^{l}\leq x\right)\\
	=&\sum_{d=1}^{\infty}\sum_{m=0}^{nt-1}\mathbf{P}\left(\frac{\tilde{Z}_{nt}^{l}}{a_{nt}}\leq x\Big|\tilde{D}_{nt}^{l,1}=d, \tilde{G}_{nt}^{l,1}=m +1\right)\cdot\mathbf{P}\left(\tilde{D}_{nt}^{l,1}=d , \tilde{G}_{nt}^{l,1}=m +1\right)\\
	=&\sum_{d =1}^{\infty}\sum_{m =0}^{nt-1}\mathbf{P}\left(\frac{\tilde{Z}_{m +n-nt,\frac{m }{m +n-nt}}^{l}}{a_{nt}}
	+\sum_{j=1}^{d}\frac{L_{m_l+n-nt,\frac{m_l}{m_l+n-nt};j}}{a_{nt}}\leq x\right)\cdot\mathbf{P}\left(\tilde{D}_{nt}^{l,1}=d , \tilde{G}_{nt}^{l,1}=m +1\right)\\
	=&\sum_{d=1}^{\infty}\sum_{m =0}^{nt-1}\mathbf{P}\left(\frac{a_{N_{nt}^{l,1},T_{nt}^{l,1}}}{a_{nt}}\frac{\tilde{Z}_{N_{nt}^{l,1},T_{nt}^{l,1}}^{l}}{a_{N_{nt}^{l,1},T_{nt}^{l,1}}}
	+\frac{b_{N_{nt}^{l,1},T_{nt}^{l,1}}}{a_{nt}}\sum_{j=1}^{\tilde{D}_{nt}^{l,1}}\frac{L_{N_{nt}^{l,1},T_{nt}^{l,1};j}}{b_{N_{nt}^{l,1},T_{nt}^{l,1}}}\leq x, \tilde{D}_{nt}^{l,1}=d , \tilde{G}_{nt}^{l,1}=m +1\right)\\
	=&\mathbf{P}\left(\frac{a_{N_{nt}^{l,1},T_{nt}^{l,1}}}{a_{nt}}\dot{Z}_{N_{nt}^{l,1},T_{nt}^{l,1}}^{l}
	+\frac{b_{N_{nt}^{l,1},T_{nt}^{l,1}}}{a_{nt}}\sum_{j=1}^{\tilde{D}_{nt}^{l,1}}\dot{L}_{N_{nt}^{l,1},T_{nt}^{l,1};j}\leq x\right).
	\end{align*}
	
	For $k>1$, the result follows by a straightforward adaptation after partitioning via the event $\bigcap_{i=1}^{k}\{\tilde{D}_{nt}^{l,i}=d_{l,i},\tilde{G}_{nt}^{l,i}=m_{l,i}+1\}$, and we omit the details.
\end{proof}

\begin{proof}[Proof of Theorem \ref{th:D-l}]
For any fixed $k\geq1$, by Theorem \ref{th:MRCA-l}, as $n\rightarrow\infty$,
\[\frac{a_{N_{nt}^{l,k},T_{nt}^{l,k}}}{a_{n,t}}\sim \frac{b_{N_{nt}^{l,k},T_{nt}^{l,k}}}{a_{n,t}}=g_{t}\left(\frac{\tilde{G}_{nt}^{l,k}}{nt}\right)\overset{\text{d}}{\rightarrow}U_{k},\]
where $\{U_{k}\}_{k\geq1}$ is a sequence of nested uniform random variables on $[0,1]$. Then by Proposition \ref{prop:nurv}, $\mathbf{E}U_k=1/2^k$.

Taking expectations on either side of \eqref{eq:dotZ-disl} yields
\[1=\mathbf{E}\left(\frac{a_{N_{nt}^{l,k},T_{nt}^{l,k}}}{a_{n,t}}\right)+\sum_{i=1}^{k}\mathbf{E}\left(\frac{b_{N_{nt}^{l,i},T_{nt}^{l,i}}}{a_{n,t}}\tilde{D}_{nt}^{l,i}\right).\]

Since $\tilde{D}_{nt}^{l,i}\geq1$ unless $\tilde{G}_{nt}^{l,i}=0$, we have, as $n\rightarrow\infty$, $\tilde{D}_{nt}^{l,i}\overset{\text{d}}{\rightarrow}1$ for all $1\leq i\leq k$. Due to the arbitrariness of $k$, $\tilde{D}_{nt}^{l,k}\overset{\text{d}}{\rightarrow}1$ for all $k\geq1$.
\end{proof}

\subsection{Limiting Distribution}
Geiger \cite{Geiger00} characterizes the exponential distribution by the following equation: if $X_{1}$ and $X_{2}$ are independent copies of a random variable $X$ with positive finite variance, and $U$ is a uniform random variable on $[0,1]$ independent of the $X_{i}$, then $X$ is exponentially distributed if and only if 
\begin{equation}\label{eq-exp}
X\overset{\text{d}}{=}U(X_{1}+X_{2}).	
\end{equation}

Since Theorem \ref{th:MRCA-l} reveals that the most recent common ancestor of the left subtree at generation $[nt]$ is uniformly distributed, and Theorem \ref{th:D-l} shows that asymptotically only one left-sibling who still have surviving descendants in  generation $[nt]$ will be produced, we can derive an expression similar to \eqref{eq-exp} for the standardized random variables, and then the convergence follows.

 To do this, recall that the $L_{2}$--Wasserstein distance between two probability measures $\mu$ and $\nu$ on $\mathbb{R}$ with finite second moment is defined as 
\begin{equation}\label{eq-W}
	d_{2}(\mu,\nu):=\inf_{X\overset{\text{d}}{=}\mu, Y\overset{\text{d}}{=}\nu}\sqrt{\mathbf{E}(X-Y)^2},
\end{equation}
where the infimum is over all pairs of random variables $X$ and $Y$, where $X$ has law $\nu$ and $Y$ has law $\mu$. And the distance has following basic properties (see Section 8 of \cite{Bickel} for detail):
\begin{equation}\label{eq-W1}
	\text{The infimum in \eqref{eq-W} is attained.}
\end{equation}
\begin{equation}\label{eq-W2}
	d_{2}(\nu_{n},\nu)\rightarrow0\Leftrightarrow\nu_{n}\rightarrow\nu \text{ weakly and } \int x^2\nu_{n}(dx)\rightarrow\int x^2\nu(dx).
\end{equation}

\begin{proof}[Proof of Theorem \ref{th:yaglom-l}]
Let $\bar{Z}_{nt}^{l}, \bar{G}_{nt}^{l}$  denote random variables with distributions
\begin{equation}\label{eq:barZ-defl}
	\mathscr{L}(\bar{Z}_{nt}^{l})=\mathscr{L}\left(\frac{\tilde{Z}_{nt}^{l}}{a_{nt}}\Big|\tilde{D}_{nt}^{l}=1\right) \quad \text{and} \quad
	\mathscr{L}(\bar{G}_{nt}^{l})=\mathscr{L}\left(\tilde{G}_{nt}^{l}\big|\tilde{D}_{nt}^{l}=1\right),
\end{equation}
where for the simplicity of notation we drop all the $1$ in the superscript $\{l,1\}$ from now on. Then using a same method as in Theorem \ref{th:decom-3l}, we have, for $n\geq1$, $0<t<1$,
\begin{equation}\label{eq:barZ-disl}
	\bar{Z}_{nt}^{l}\overset{\text{d}}{=}\frac{a_{\bar{N}_{nt}^{l},\bar{T}_{nt}^{l}}}{a_{n,t}}\dot{Z}_{\bar{N}_{nt}^{l},\bar{T}_{nt}^{l}}^{l}
	+\frac{b_{\bar{N}_{nt}^{l},\bar{T}_{nt}^{l}}}{a_{n,t}}\dot{L}_{\bar{N}_{nt}^{l},\bar{T}_{nt}^{l}},
\end{equation}
where $\bar{N}_{nt}^{l}:=\bar{G}_{nt}^{l}-1+n-nt$, $\bar{T}_{nt}^{l}:=(\bar{G}_{nt}^{l}-1)/(\bar{G}_{nt}^{l}-1+n-nt)$. And for all $i\geq0$, $s\in(0,1)$, $\dot{Z}^{l}_{i,s}$, $\dot{L}_{i,s}$ are as in Theorem \ref{th:decom-3l} and independent of $\bar{G}_{nt}^{l}$. 

Let $X$ be a random variable with exponential distribution with mean $1$, we will show that the $d_{2}$--distance between $\mathscr{L}(\dot{Z}_{nt}^{l})$ and $\mathscr{L}(X)$ 
\begin{equation}\label{eq4}
		d_{2}(\mathscr{L}(\dot{Z}_{nt}^{l}),\mathscr{L}(X))^{2}\leq
		d_{2}(\mathscr{L}(\dot{Z}_{nt}^{l}),\mathscr{L}(\bar{Z}_{nt}^{l}))^{2}
		+d_{2}(\mathscr{L}(\bar{Z}_{nt}^{l}),\mathscr{L}(X))^{2}
\end{equation}
goes to $0$ as $n\rightarrow\infty$, where $\dot{Z}_{nt}^{l}$ and $\bar{Z}_{nt}^{l}$ are defined in \eqref{eq:dotZ-defl}, \eqref{eq:barZ-defl}, and have distributions \eqref{eq:dotZ-disl}, \eqref{eq:barZ-disl}.

$(i)$ The $d_{2}$--distance between $\dot{Z}_{nt}^{l}$ and $\bar{Z}_{nt}^{l}$:

Since by Theorem \ref{th:moment-l2} and  Lemma \ref{lem:Lnt}, for all $0<t<1$, as $n\rightarrow\infty$,
\[\mathbf{E}(\dot{Z}_{nt}^{l})^{2}=a_{nt}^{-2}\mathbf{E}(\tilde{Z}_{nt}^{l})^{2}\rightarrow2,\quad \mathbf{E}(\dot{L}_{nt})^{2}
		=b_{nt}^{-2}\mathbf{E}(L_{nt})^{2}\rightarrow2,\]
we have,
{\small
\begin{align*}
		\mathbf{E}(\bar{Z}_{nt}^{l})^{2}
		=&\mathbf{E}\left(\frac{a_{\bar{N}_{nt}^{l}\bar{T}_{nt}^{l}}}{a_{nt}}\right)^{2}
		\mathbf{E}\left(\dot{Z}_{\bar{N}_{nt}^{l}\bar{T}_{nt}^{l}}^{l}\right)^{2}
		+\mathbf{E}\left(\frac{b_{\bar{N}_{nt}^{l}\bar{T}_{nt}^{l}}}{a_{nt}}\right)^{2}
		\mathbf{E}\left(\dot{L}_{\bar{N}_{nt}^{l}\bar{T}_{nt}^{l}}\right)^{2}
		+\mathbf{E}\left(\frac{2a_{\bar{N}_{nt}^{l}\bar{T}_{nt}^{l}}b_{\bar{N}_{nt}^{l}\bar{T}_{nt}^{l}}}{a_{nt}^{2}}\right)\\
	\sim& \mathbf{E}U^{2}\cdot 
		\left[\mathbf{E}\left(\dot{Z}_{\bar{N}_{nt}^{l}\bar{T}_{nt}^{l}}^{l}\right)^{2}
		+\mathbf{E}\left(\dot{L}_{\bar{N}_{nt}^{l}\bar{T}_{nt}^{l}}\right)^{2}+2\right]\rightarrow2.
\end{align*}}

Consequently, by Theorem \ref{th:D-l} and \eqref{eq-W2}, 
\begin{equation}\label{eq4-0}
	d_{2}(\mathscr{L}(\dot{Z}_{nt}^{l}),\mathscr{L}(\bar{Z}_{nt}^{l}))^{2}\rightarrow0.
\end{equation}

$(ii)$ The $d_{2}$--distance between $\bar{Z}_{nt}^{l}$ and $X$:

Let $X_{1}$, $X_{2}$ and $U$ be independent random variables where $X_{i}$ has exponential distribution with mean $1$ and $U$ is uniform on $[0,1]$, then by \eqref{eq-exp}, $X=U(X_{1}+X_{2})$ is also exponential distributed with mean $1$. 

By \eqref{eq-W1} we can choose versions $\{\tilde{G}_{nt}^{l},\dot{Z}_{nt}^{l}, \dot{L}_{nt}\}_{n\geq0,t\in(0,1)}$, such that
\begin{align}
	\mathbf{E}\left(\frac{n^{-1}\tilde{G}_{nt}^{l}}{t(n^{-1}\tilde{G}_{nt}^{l}+1-t)}-U\right)^{2}&=d_{2}\left(\mathscr{L}\left(\frac{n^{-1}\tilde{G}_{nt}^{l}}{t(n^{-1}\tilde{G}_{nt}^{l}+1-t)}\right),\mathscr{L}(U)\right)^{2},\label{eq4-1}\\
	\mathbf{E}\left(\dot{Z}_{nt}^{l}-X_{1}\right)^{2}&=d_{2}\left(\mathscr{L}(\dot{Z}_{nt}^{l}),\mathscr{L}(X)\right)^{2},\label{eq4-2}\\
	\mathbf{E}\left(\dot{L}_{nt}-X_{2}\right)^{2}&=d_{2}\left(\mathscr{L}(\dot{L}_{nt}),\mathscr{L}(X)\right)^{2},\label{eq4-3}
\end{align}
where for any fixed $0<t<1$, the right hand in \eqref{eq4-1} and \eqref{eq4-3} tend to 0 in view of  Theorem \ref{th:MRCA-l}, Lemma \ref{lem:Lnt}, \eqref{eq-W2}, Skorokhod representation theorem, $n^{-1}\tilde{G}_{nt}^{l}\leq t$ and dominated convergence theorem.  Then
	\begin{align*}
		d_{2}\left(\mathscr{L}(\bar{Z}_{nt}^{l}),\mathscr{L}(X)\right)^{2}
		\leq &\mathbf{E}\left(\bar{Z}_{nt}^{l}-X\right)^{2}\\
		=&\mathbf{E}\left(\frac{a_{\bar{N}_{nt}^{l}\bar{T}_{nt}^{l}}}{a_{nt}}\dot{Z}_{\bar{N}_{nt}^{l}\bar{T}_{nt}^{l}}^{l}-UX_{1}+\frac{b_{\bar{N}_{nt}^{l}\bar{T}_{nt}^{l}}}{a_{nt}}\dot{L}_{\bar{N}_{nt}^{l}\bar{T}_{nt}^{l}}-UX_{2}\right)^{2}\\
		=&\mathbf{E}\left(\frac{a_{\bar{N}_{nt}^{l}\bar{T}_{nt}^{l}}}{a_{nt}}\dot{Z}_{\bar{N}_{nt}^{l}\bar{T}_{nt}^{l}}^{l}-UX_{1}\right)^{2}
		+\mathbf{E}\left(\frac{b_{\bar{N}_{nt}^{l}\bar{T}_{nt}^{l}}}{a_{nt}}\dot{L}_{\bar{N}_{nt}^{l}\bar{T}_{nt}^{l}}-UX_{2}\right)^{2}\\
		&+2\mathbf{E}\left(\frac{a_{\bar{N}_{nt}^{l}\bar{T}_{nt}^{l}}}{a_{nt}}\dot{Z}_{\bar{N}_{nt}^{l}\bar{T}_{nt}^{l}}^{l}-UX_{1}\right)\cdot\left(\frac{b_{\bar{N}_{nt}^{l}\bar{T}_{nt}^{l}}}{a_{nt}}\dot{L}_{\bar{N}_{nt}^{l}\bar{T}_{nt}^{l}}-UX_{2}\right)\\
		:=&I_{1}+I_{2}+2I_3.
	\end{align*}
	
For the third part, since $\mathbf{E}X_{1}=\mathbf{E}X_{2}=\mathbf{E}\dot{Z}_{\bar{N}_{nt}^{l}\bar{T}_{nt}^{l}}^{l}=\mathbf{E}\dot{L}_{\bar{N}_{nt}^{l}\bar{T}_{nt}^{l}}=1$, by \eqref{eq4-1}, as $n\rightarrow\infty$,
	\begin{align*}
		I_{3}
		=&\mathbf{E}\left(\frac{a_{\bar{N}_{nt}^{l}\bar{T}_{nt}^{l}}b_{\bar{N}_{nt}^{l}\bar{T}_{nt}^{l}}}{a_{nt}^{2}}\dot{Z}_{\bar{N}_{nt}^{l}\bar{T}_{nt}^{l}}^{l}\dot{L}_{\bar{N}_{nt}^{l}\bar{T}_{nt}^{l}}\right)
		-\mathbf{E}\left(\frac{a_{\bar{N}_{nt}^{l}\bar{T}_{nt}^{l}}}{a_{nt}}\dot{Z}_{\bar{N}_{nt}^{l}\bar{T}_{nt}^{l}}^{l}UX_{2}\right)\\
		&-\mathbf{E}\left(UX_{1}\frac{b_{\bar{N}_{nt}^{l}\bar{T}_{nt}^{l}}}{a_{nt}}\dot{L}_{\bar{N}_{nt}^{l}\bar{T}_{nt}^{l}}\right)
		+\mathbf{E}\left(U^{2}X_{1}X_{2}\right)\\
		=&\mathbf{E}\left(\frac{a_{\bar{N}_{nt}^{l}\bar{T}_{nt}^{l}}b_{\bar{N}_{nt}^{l}\bar{T}_{nt}^{l}}}{a_{nt}^{2}}\right)
		-\mathbf{E}\left(\frac{a_{\bar{N}_{nt}^{l}\bar{T}_{nt}^{l}}}{a_{nt}}U\right)
		-\mathbf{E}\left(U\frac{b_{\bar{N}_{nt}^{l}\bar{T}_{nt}^{l}}}{a_{nt}}\right)
		+\mathbf{E}\left(U^{2}\right)\\
		=&\mathbf{E}\left(\frac{a_{\bar{N}_{nt}^{l}\bar{T}_{nt}^{l}}}{a_{nt}}-U\right)\left(\frac{b_{\bar{N}_{nt}^{l}\bar{T}_{nt}^{l}}}{a_{nt}}-U\right)\\
		=&\mathbf{E}\left(\frac{n^{-1}\tilde{G}_{nt}^{l}}{t(n^{-1}\tilde{G}_{nt}^{l}+1-t)}-U\right)^{2}+o(1)\\
		\rightarrow & \,0.
	\end{align*}

For the second part, by \eqref{eq4-3} and dominated convergence theorem, as $n\rightarrow\infty$,
	\begin{align*}
		I_{2}
		=&\mathbf{E}\left(\frac{b_{\bar{N}_{nt}^{l}\bar{T}_{nt}^{l}}}{a_{nt}}\left(\dot{L}_{\bar{N}_{nt}^{l}\bar{T}_{nt}^{l}}-X_{2}\right)+\left(\frac{b_{\bar{N}_{nt}^{l}\bar{T}_{nt}^{l}}}{a_{nt}}-U\right)X_{2}\right)^{2}\\
		=&\mathbf{E}\left(\frac{b_{\bar{N}_{nt}^{l}\bar{T}_{nt}^{l}}}{a_{nt}}\left(\dot{L}_{\bar{N}_{nt}^{l}\bar{T}_{nt}^{l}}-X_{2}\right)\right)^{2}+\mathbf{E}\left(\frac{b_{\bar{N}_{nt}^{l}\bar{T}_{nt}^{l}}}{a_{nt}}-U\right)^{2}\cdot \mathbf{E}\left(X_{2}\right)^{2}\\
		&+2\left(1-\mathbf{E}(X_{2})^{2}\right)\cdot \mathbf{E}\left(\frac{b_{\bar{N}_{nt}^{l}\bar{T}_{nt}^{l}}}{a_{nt}}\left(\frac{b_{\bar{N}_{nt}^{l}\bar{T}_{nt}^{l}}}{a_{nt}}-U\right)\right)\\
		=&\mathbf{E}\left(\frac{b_{\bar{N}_{nt}^{l}\bar{T}_{nt}^{l}}}{a_{nt}}\left(\dot{L}_{\bar{N}_{nt}^{l}\bar{T}_{nt}^{l}}-X_{2}\right)\right)^{2}+o(1)\\
		\rightarrow&0.
	\end{align*}
	
For the first part, similarly,
\[
	I_{1}=\int_{0}^{t}\left(\frac{a_{n(s+1-t),s/(s+1-t)}}{a_{nt}}\right)\mathbf{E}\left(\dot{Z}_{n(s+1-t),s/(s+1-t)}^{l}-X_{1}\right)^{2}\mathbf{P}\left(\frac{\tilde{G}_{nt}^{l}}{n}\in\text{d}s\right)+o(1).\\
\]

$(iii)$ The $d_{2}$--distance between $\dot{Z}_{nt}^{l}$ and $X$:
 
By combining the above equations with  \eqref{eq4}, \eqref{eq4-0} and \eqref{eq4-2}, we have
\begin{align*}
	&\sup_{0<t<1}\varlimsup_{n\rightarrow\infty}d_{2}\left(\mathscr{L}(\dot{Z}_{nt}^{l}),\mathscr{L}(X)\right)^{2}\\
	\leq&\sup_{0<t<1}\varlimsup_{n\rightarrow\infty}d_{2}\left(\mathscr{L}(\bar{Z}_{nt}^{l}),\mathscr{L}(X)\right)^{2}\\
	=&\sup_{0<t<1}\varlimsup_{n\rightarrow\infty}\int_{0}^{t}\left(\frac{a_{n(s+1-t),\frac{s}{s+1-t}}}{a_{nt}}\right)\mathbf{E}\left(\dot{Z}_{n(s+1-t),\frac{s}{s+1-t}}^{l}-X_{1}\right)^{2}\mathbf{P}\left(\frac{\tilde{G}_{nt}^{l}}{n}\in\text{d}s\right)\\
	\leq&\sup_{0<t<1}\int_{0}^{t}\varlimsup_{n\rightarrow\infty}\left(\frac{a_{n(s+1-t),\frac{s}{s+1-t}}}{a_{nt}}\right)\varlimsup_{n\rightarrow\infty}d_{2}\left(\mathscr{L}(\dot{Z}_{n(s+1-t),\frac{s}{s+1-t}}^{l}),\mathscr{L}(X)\right)^2\mathbf{P}\left(\frac{\tilde{G}_{nt}^{l}}{n}\in\text{d}s\right)\\
	\leq&\sup_{0<t<1}\varlimsup_{n\rightarrow\infty}d_{2}\left(\mathscr{L}(\dot{Z}_{nt}^{l}),\mathscr{L}(X)\right)^2\cdot\mathbf{E}U^2.
\end{align*}

Since $\mathbf{E}U^2=1/3$ and $\limsup_{n\rightarrow\infty}d_{2}\left(\mathscr{L}(\dot{Z}_{nt}^{l}),\mathscr{L}(X)\right)<\infty$, we have, for any fixed $0<t<1$, 
\[\lim_{n\rightarrow\infty}d_{2}\left(\mathscr{L}(\dot{Z}_{nt}^{l}),\mathscr{L}(X)\right)^{2}=0,\]
and the theorem follows.
\end{proof}

\section{Genealogy of the Right Subtree}\label{section:3}
The branching mechanism of the right subtree follows the ordinary Galton–Watson process. Consequently, the analysis is similar to (and  simpler than) that of the left part; we therefore omit the details and present only the key steps.
\subsection{Moments}
By Theorem \ref{th:decom-1},
\begin{equation*}
	\tilde{Z}_{nt}^{r}\overset{\text{d}}{=}1+\sum_{j=1}^{X_{n-nt+1}}Z_{nt,j}^{r,0}+\sum_{j=1}^{X_{n-nt+2}}Z_{nt,j}^{r,1}+\cdots+\sum_{j=1}^{X_{n}}Z_{nt,j}^{r,nt-1},
\end{equation*}
where $Z_{nt;j}^{r,i}\overset{\text{d}}{=}Z_{i}$. Since $\mathbf{E}Z_{i}=1$ and $\mathbf{E}Z_{i}^{2}=1+i\sigma^2$, we have,
\begin{theorem}\label{th:moment-r}
	Suppose $\alpha=1$ and $\sigma^{2}<\infty$, then  for any fixed $0<t<1$,
	\[
		\lim_{n\rightarrow\infty}n^{-1}\mathbf{E}{\tilde{Z}_{nt}^{r}}=\frac{t\sigma^{2}}{2},\quad
		\lim_{n\rightarrow\infty}n^{-2}\mathbf{E}(\tilde{Z}_{nt}^{r})^{2}=\frac{(t\sigma^{2})^{2}}{2}.
	\]
\end{theorem}
\subsection{Coalescent Times}
\begin{proof}[Proof of Theorem \ref{th:MRCA-r}]
Let $\{U_{k}\}_{k\geq1}$ be a nested uniform random variable sequence on $[0,1]$ defined in Definition \ref{def:nurv}. Since $(\tilde{Z}_{nt}^{r,i})_{i\geq0}$ have non-negative independent increments, by \eqref{eq3-3},
\[\mathbf{P}\left(\tilde{G}_{nt}^{r,1}\leq nx\right)=\prod_{i=nx}^{nt-1}\mathbf{P}\left(\tilde{Z}_{nt}^{r,i}=\tilde{Z}_{nt}^{r,i+1}\right)
	\sim\prod_{i=nx}^{nt-1}\left(1-\frac{1}{i}\right)\rightarrow\frac{x}{t}
\]
for $0\leq x\leq t$. Thus as $n\rightarrow\infty$, $\tilde{G}_{nt}^{r,1}/nt\overset{\text{d}}{\rightarrow}U_{1}$.

Assume that for some $k\geq1$, $\tilde{G}_{nt}^{r,k}/nt\overset{\text{d}}{\rightarrow}U_{k}$. Then for $0\leq x\leq t$, if $ny\leq nx$,
\[\mathbf{P}\left(\tilde{G}_{nt}^{r,k+1}\leq nx\Big|\tilde{G}_{nt}^{r,k}=ny\right)=1;\]
if $ny>nx$,  
\[\mathbf{P}\left(\tilde{G}_{nt}^{r,k+1}\leq nx\Big|\tilde{G}_{nt}^{r,k}=ny\right)=\mathbf{P}\left(\tilde{Z}_{nt}^{r,nx}=\cdots=\tilde{Z}_{nt}^{r,ny-1}\right)\sim \prod_{i=nx}^{ny-2}\left(1-\frac{1}{i}\right)\rightarrow\frac{x}{y}.\]
Thus as $n\rightarrow\infty$, $\tilde{G}_{nt}^{r,k+1}/nt\overset{\text{d}}{\rightarrow}U[0,U_{k}]=U_{k+1}$, and the theorem follows by induction.
\end{proof}
\begin{corollary}
	Suppose $\alpha=1$ and $\sigma^{2}<\infty$, then for any fixed $0<t<1$, $k\in\mathbb{N}_{+}$, and $0\leq x\leq t$,
\[\lim_{n\rightarrow\infty}\mathbf{P}\left(\frac{\tilde{G}_{nt}^{r,k}}{n}\leq x\right)=\frac{x}{t}\sum_{m=0}^{k-1}\frac{1}{m!}\left(\ln\frac{t}{x}\right)^{m}.
	\]
\end{corollary}

\subsection{Siblings of the Spine}
Let $R_n$ be the random variable with distribution
\[\mathscr{L}(R_{n})=\mathscr{L}(Z_{n}|Z_{n}>0),\]
the the limiting distribution is given by Yaglom's Theorem. 

Define
\[
	\dot{Z}_{nt}^{r}:=\frac{\tilde{Z}_{nt}^{r}}{c_{nt}} \quad\text{and}\quad \dot{R}_{n}:=\frac{R_{n}}{d_{n}},
\]
where $c_{nt}:=\mathbf{E}\tilde{Z}_{nt}^{r}$, $d_{n}:=\mathbf{E}R_{n}$ are the mean values. Then  $=\mathbf{E}{\dot{Z}_{nt}^{r}}=\mathbf{E}{\dot{R}_{n}}=1$, and by Theorem \ref{th:moment-r}, $\mathbf{P}(Z_{n}>0)\sim 2/(\sigma^{2}n)$, for any fixed $0<t<1$,
\[\lim_{n\rightarrow\infty}\frac{c_{nt}}{nt}=\lim_{n\rightarrow\infty}\frac{d_{n}}{n}=\frac{\sigma^{2}}{2}.\]

Then we have the following expression for the distribution of the standardized random variable, from which the limiting distribution of $\tilde{D}_{nt}^{r,k}$ is derived by taking expectations.
\begin{theorem}\label{th:decom-3r}
	For all $n\geq0$ and $0<t<1$, $k\geq1$,
\begin{equation}\label{eq:dot-disr}
	\dot{Z}_{nt}^{r}\overset{\text{d}}{=}\frac{c_{N_{nt}^{r,k}T_{nt}^{r,k}}}{c_{nt}}\dot{Z}_{N_{nt}^{r,k}T_{nt}^{r,k}}^{r}+\sum_{i=1}^{k}\frac{d_{\tilde{G}_{nt}^{r,k}-1}}{c_{nt}}\sum_{j=1}^{\tilde{D}_{nt}^{r,i}}\dot{R}_{\tilde{G}_{nt}^{r,k}-1;j},
\end{equation}
where $N_{nt}^{r,k}:=\tilde{G}_{nt}^{r,k}-1+n-nt$, $T_{nt}^{r,k}:=(\tilde{G}_{nt}^{r,k}-1)/(\tilde{G}_{nt}^{r,k}-1+n-nt)$. And for all $m\geq0$, $0<s<1$, $\{\dot{R}_{m;j}\}_{j\geq1}$ are independent copies of $\dot{R}_{m}$; $\{\dot{R}_{ms;j},\dot{Z}_{m,s}^{r}\}_{m,s,j}$ are independent and independent of $\{\tilde{D}_{nt}^{r,k},\tilde{G}_{nt}^{r,k}\}_{k}$.
\end{theorem}

\begin{proof}[Proof of Theorem \ref{th:D-r}]
By Theorem  \ref{th:MRCA-r}, for all $k\geq1$, as $n\rightarrow\infty$,
\[\frac{c_{N_{nt}^{r,k},T_{nt}^{r,k}}}{c_{n,t}}\sim \frac{d_{\tilde{G}_{nt}^{r,k}-1}}{c_{n,t}}\sim\frac{\tilde{G}_{nt}^{r,k}}{nt}\overset{\text{d}}{\rightarrow}U_{k},\]
where $\{U_{k}\}_{k\geq1}$ is a sequence of nested uniform random variables on $[0,1]$. Then the the theorem follows by taking expectations on either side of the equation \eqref{eq:dot-disr}.
\end{proof}

\subsection{Limiting Distribution}

\begin{proof}[Proof of Theorem \ref{th:yaglom-r}]
Let $\bar{Z}_{nt}^{r}, \bar{G}_{nt}^{r}$  denote random variables with distributions
\[\mathscr{L}(\bar{Z}_{nt}^{r})=\mathscr{L}\left(\frac{\tilde{Z}_{nt}^{r}}{c_{nt}}\Big|\tilde{D}_{nt}^{r}=1\right)\quad\text{and}\quad \mathscr{L}(\bar{G}_{nt}^{r})=\mathscr{L}\left(\tilde{G}_{nt}^{r}\big|\tilde{D}_{nt}^{r}=1\right),\]
then we have, for $n\geq1$ and $0<t<1$,
\[\bar{Z}_{nt}^{r}\overset{\text{d}}{=}\frac{c_{\bar{N}_{nt}^{r}\bar{T}_{nt}^{r}}}{c_{nt}}\dot{Z}_{\bar{N}_{nt}^{r}\bar{T}_{nt}^{r}}^{r}
	+\frac{d_{\bar{G}_{nt}^{r}-1}}{c_{nt}}\dot{R}_{\bar{G}_{nt}^{r}-1},\]
where $\bar{N}_{nt}^{r}:=\bar{G}_{nt}^{r}-1+n-nt$, $\bar{T}_{nt}^{r}:=(\bar{G}_{nt}^{r}-1)/(\bar{G}_{nt}^{r}-1+n-nt)$. And for all $i\geq0$, $s\in(0,1)$, $\dot{Z}^{r}_{r,s}$, $\dot{R}_{is}$ are as in Theorem \ref{th:decom-3r} and independent of $\bar{G}_{nt}^{r}$. 

For the convergence, just notice that by Theorem \ref{th:moment-r} and $\mathbf{P}(Z_{n}>0)\sim 2/(\sigma^{2}n)$, 
\[\mathbf{E}(\dot{Z}_{nt}^{r})^{2}=c_{nt}^{-2}\mathbf{E}(\tilde{Z}_{nt}^{r})^{2}\rightarrow2,\quad \mathbf{E}(\dot{R}_{n})^{2}=\frac{\mathbf{E}(Z_{n}^{2}|Z_{n}>0)}{d_{n}^{2}}=\frac{1+\sigma^{2}n}{d_{n}^{2}\mathbf{P}(Z_{n}>0)}\rightarrow2,\]
then $\mathbf{E}(\bar{Z}_{nt}^{r})^{2}\rightarrow2$, and $d_{2}(\mathscr{L}(\dot{Z}_{nt}^{r}),\mathscr{L}(\bar{Z}_{nt}^{r}))^{2}\rightarrow0$. And the proof of is $d_{2}(\mathscr{L}(\dot{Z}_{nt}^{r}),\mathscr{L}(X))^{2}\rightarrow0$ analogous to the proof of Theorem \ref{th:yaglom-l}.

\end{proof}

\section{Asymptotic Independence of the Two Subtrees}\label{section:4}

In  this section, our main objective is to prove the next theorem for the joint distribution of coalescent time sequences $\{\tilde{G}_{nt}^{l,k}\}_{k}$ and $\{\tilde{G}_{nt}^{r,k}\}_{k}$, by which Theorem \ref{th:MRCA-in} is a direct result.
\begin{theorem}
	Suppose $\alpha=1$ and $\sigma^{2}<\infty$, then for any fixed $0<t<1$, $k_{l},k_{r}\in\mathbb{N}_{+}$, and $x,y\in[0,t]$,
	{\small
	\[
	\lim_{n\rightarrow\infty}\mathbf{P}(\tilde{G}_{nt}^{l,k_{l}}\leq nx, \tilde{G}_{nt}^{r,k_{r}}\leq ny)
	=\frac{xy}{t^{2}(x+1-t)}
	\left(\sum_{w=0}^{k_{l}-1}\frac{1}{w!}\left(\ln\frac{t(x+1-t)}{x}\right)^{w}\right)
	\left(\sum_{m=0}^{k_{r}-1}\frac{1}{m!}\left(\ln\frac{t}{y}\right)^{m}\right).
	\]}
\end{theorem}
\begin{proof}
	Notice that for all $k_{l},k_{r}\in\mathbb{N}_{+}$ and $x,y\in[0,t]$, the event $\{\tilde{G}_{nt}^{l,k_{l}}\leq nx, \tilde{G}_{nt}^{r,k_{r}}\leq ny\}$ is equivalent to 
\begin{equation*}
	\left\{\big|\{nx\leq i<nt:\tilde{Z}_{nt}^{l,i+1}>\tilde{Z}_{nt}^{l,i}\}\big|<k_{l}\right\}\cap\left\{\big|\{ny\leq i<nt:\tilde{Z}_{nt}^{r,i+1}>\tilde{Z}_{nt}^{r,i}\}\big|<k_{r}\right\},
\end{equation*}
thus
\[\mathbf{P}\left(\tilde{G}_{nt}^{l,k_l}\leq nx, \tilde{G}_{nt}^{r,k_{r}}\leq ny\right)=P_{n,1}+P_{n,2}+P_{n,3}+P_{n,4},\]
where
\begin{align*}
	&P_{n,1}:=\mathbf{P}\left(\big|\{nx\leq i<nt:\tilde{Z}_{nt}^{l,i+1}>\tilde{Z}_{nt}^{l,i}\}\big|=0,
	\big|\{ny\leq i<nt:\tilde{Z}_{nt}^{r,i+1}>\tilde{Z}_{nt}^{r,i}\}\big|=0\right),\\
	&P_{n,2}:=\sum_{m=1}^{k_{r}-1}\mathbf{P}\left(\big|\{nx\leq i<nt:\tilde{Z}_{nt}^{l,i+1}>\tilde{Z}_{nt}^{l,i}\}\big|=0,
	\big|\{ny\leq i<nt:\tilde{Z}_{nt}^{r,i+1}>\tilde{Z}_{nt}^{r,i}\}\big|=m\right),\\
	&P_{n,3}:=\sum_{w=1}^{k_{l}-1}\mathbf{P}\left(\big|\{nx\leq i<nt:\tilde{Z}_{nt}^{l,i+1}>\tilde{Z}_{nt}^{l,i}\}\big|=w,
	\big|\{ny\leq i<nt:\tilde{Z}_{nt}^{r,i+1}>\tilde{Z}_{nt}^{r,i}\}\big|=0\right),\\
	&P_{n,4}:=\sum_{w=1}^{k_{l}-1}\sum_{m=1}^{k_{r}-1}\mathbf{P}\left(\big|\{nx\leq i<nt:\tilde{Z}_{nt}^{l,i+1}>\tilde{Z}_{nt}^{l,i}\}\big|=w,
	\big|\{ny\leq i<nt:\tilde{Z}_{nt}^{r,i+1}>\tilde{Z}_{nt}^{r,i}\}\big|=m\right),
\end{align*}
and we will estimate the limits of these four parts respectively.

$(i)$ The limit of $P_{n,1}$:

When $0\leq x=y\leq t$, since $(\tilde{Z}_{nt}^{l,i})_{i\geq0}$ and $(\tilde{Z}_{nt}^{r,i})_{i\geq0}$ both have non-negative independent increments, by \eqref{eq3-5}, we have, as $n\rightarrow\infty$,
	\[
	P_{n,1}=\prod_{i=nx}^{nt-1}\mathbf{P}(\tilde{Z}_{nt}^{l,i}=\tilde{Z}_{nt}^{l,i+1}, \tilde{Z}_{nt}^{r,i}=\tilde{Z}_{nt}^{r,i+1})
		\sim\prod_{i=nx}^{nt-1}\left[1-\left(\frac{2}{i}-\frac{1}{i+n-nt}\right)\right]\rightarrow\frac{x^2}{t^2(x+1-t)}.
	\]

When $0\leq x<y\leq t$, by \eqref{eq3-1}, \eqref{eq3-5} and the proof of Theorem \ref{th:MRCA-l}, as $n\rightarrow\infty$,
\begin{align*}
	P_{n,1}
	=&\prod_{i=nx}^{ny-1}\mathbf{P}(\tilde{Z}_{nt}^{l,i}=\tilde{Z}_{nt}^{l,i+1})
		\prod_{i=ny}^{nt-1}\mathbf{P}(\tilde{Z}_{nt}^{l,i}=\tilde{Z}_{nt}^{l,i+1}, \tilde{Z}_{nt}^{r,i}=\tilde{Z}_{nt}^{r,i+1})\\
	\sim&\prod_{i=nx}^{ny-1}\left(1-\frac{1}{i}+\frac{1}{i+n-nt}\right)
		\prod_{i=ny}^{nt-1}\left(1-\frac{2}{i}+\frac{1}{i+n-nt}\right)\\
	\rightarrow&\frac{x(y+1-t)}{y(x+1-t)}\cdot\frac{y^{2}}{t^2(y+1-t)}=\frac{xy}{t^2(x+1-t)}.
\end{align*}
		
When $0\leq y<x\leq t$, by \eqref{eq3-3}, \eqref{eq3-5} and the proof of Theorem \ref{th:MRCA-r}, as $n\rightarrow\infty$,
\begin{align*}
	P_{n,1}
	=&\prod_{i=ny}^{nx-1}\mathbf{P}(\tilde{Z}_{nt}^{r,i}=\tilde{Z}_{nt}^{r,i+1})\prod_{i=nx}^{nt-1}\mathbf{P}(\tilde{Z}_{nt}^{l,i}=\tilde{Z}_{nt}^{l,i+1}, \tilde{Z}_{nt}^{r,i}=\tilde{Z}_{nt}^{r,i+1})\\
	\sim&\prod_{i=ny}^{nx-1}\left(1-\frac{1}{i}\right)
		\prod_{i=nx}^{nt-1}\left(1-\frac{2}{i}+\frac{1}{i+n-nt}\right)\\
	\rightarrow&\frac{y}{x}\cdot\frac{x^{2}}{t^2(x+1-t)}=\frac{xy}{t^2(x+1-t)}.
\end{align*}

Thus for for all $x,y\in[0,t]$,
\begin{equation}\label{eq3-n1}
	\lim_{n\rightarrow\infty}P_{n,1}=\frac{xy}{t^2(x+1-t)}.
\end{equation}

$(ii)$ The limit of $P_{n,2}$:

When $0\leq x<y\leq t$, $P_{n,2}$ can be represented as
\begin{align*}
	&
	\sum_{m=1}^{k_{r}-1}
	\sum_{ny\leq v_1<\cdots<v_m<nt}
	\mathbf{P}
	\left(
	\bigcap_{nx\leq s<nt}
	\left\{\tilde{Z}_{nt}^{l,s+1}=\tilde{Z}_{nt}^{l,s}\right\}
	\bigcap_{s\in \{v_1,\cdots,v_m\}}
	\left\{\tilde{Z}_{nt}^{r,s+1}>\tilde{Z}_{nt}^{r,s}\right\}
	\right.
	\\
	&\left.
	\qquad\qquad\qquad\qquad\qquad\quad
	\bigcap_{\substack{ny\leq s< nt\\s\notin\{v_1,\cdots,v_m\}}}
	\left\{\tilde{Z}_{nt}^{r,s+1}=\tilde{Z}_{nt}^{r,s}\right\}
	\right)\\
	=&
	\sum_{m=1}^{k_{r}-1}
	\sum_{ny\leq v_1<\cdots<v_m<nt}
	\left(
	\prod_{nx\leq s<ny}
	\mathbf{P}
	\left(\tilde{Z}_{nt}^{l,s+1}=\tilde{Z}_{nt}^{l,s}\right)
	\prod_{s\in \{v_1,\cdots,v_m\}}
	\mathbf{P}
	\left(\tilde{Z}_{nt}^{l,s+1}=\tilde{Z}_{nt}^{l,s},
	\right.
	\right.\\
	&\left.
	\left.
	\qquad\qquad\qquad\qquad\qquad\qquad
	\tilde{Z}_{nt}^{r,s+1}>\tilde{Z}_{nt}^{r,s}\right)
	\prod_{\substack{ny\leq s< nt\\s\notin \{v_1,\cdots,v_m\}}}
	\mathbf{P}
	\left(\tilde{Z}_{nt}^{l,s+1}=\tilde{Z}_{nt}^{l,s},\tilde{Z}_{nt}^{r,s+1}=\tilde{Z}_{nt}^{r,s}\right)
	\right)
	\\
	=&
	\sum_{m=1}^{k_{r}-1}
	\sum_{ny\leq v_1<\cdots<v_m<nt}
	\left(
	\prod_{s=nx}^{ny-1}
	\mathbf{P}
	\left(\tilde{Z}_{nt}^{l,s+1}=\tilde{Z}_{nt}^{l,s}\right)
	\prod_{s\in \{v_1,\cdots,v_m\}}
	\frac{
	\mathbf{P}
	\left(\tilde{Z}_{nt}^{l,s+1}=\tilde{Z}_{nt}^{l,s}, \tilde{Z}_{nt}^{r,s+1}>\tilde{Z}_{nt}^{r,s}\right)}
	{\mathbf{P}
	\left(\tilde{Z}_{nt}^{l,s+1}=\tilde{Z}_{nt}^{l,s}, \tilde{Z}_{nt}^{r,s+1}=\tilde{Z}_{nt}^{r,s}\right)}
	\right.
	\\
	&\left.
	\qquad\qquad\qquad\qquad\qquad\quad
	\prod_{s=ny}^{nt-1}
	\mathbf{P}
	\left(\tilde{Z}_{nt}^{l,s+1}=\tilde{Z}_{nt}^{l,s},\tilde{Z}_{nt}^{r,s+1}=\tilde{Z}_{nt}^{r,s}\right)
	\right)
	\\
	=&
	P_{n,1}
	\cdot
	\sum_{m=1}^{k_{r}-1}
	\sum_{ny\leq v_1<\cdots<v_m<nt}
	\prod_{s\in \{v_1,\cdots,v_m\}}
	\frac{
	\mathbf{P}
	\left(\tilde{Z}_{nt}^{l,s+1}=\tilde{Z}_{nt}^{l,s}, \tilde{Z}_{nt}^{r,s+1}>\tilde{Z}_{nt}^{r,s}\right)}
	{\mathbf{P}
	\left(\tilde{Z}_{nt}^{l,s+1}=\tilde{Z}_{nt}^{l,s}, \tilde{Z}_{nt}^{r,s+1}=\tilde{Z}_{nt}^{r,s}\right)},
\end{align*}
thus by \eqref{eq3-4} and \eqref{eq3-5},
\begin{align*}
	\lim_{n\rightarrow\infty}P_{n,2}
	=&\lim_{n\rightarrow\infty}P_{n,1}\cdot\lim_{n\rightarrow\infty} \sum_{m=1}^{k_{r}-1}\sum_{ny\leq v_{1}<\cdots<v_{m}< nt}\left(\frac{1}{v_1v_2\cdots v_m}\right)\\
	=&\frac{xy}{t^2(x+1-t)}\cdot\sum_{m=1}^{k_{r}-1}\left(\int_{y\leq s_{1}<\cdots<s_{m}< t}\frac{1}{s_{1}\cdots s_{m}}\text{d}s_{1}\cdots \text{d}s_{m}\right)\\
	=&\frac{xy}{t^2(x+1-t)}\cdot\sum_{m=1}^{k_{r}-1}\frac{1}{m!}\left(\int_{y}^{t}\frac{1}{s}\text{d}s\right)^{m}\\
	=&\frac{xy}{t^{2}(x+1-t)}\cdot \sum_{m=1}^{k_{r}-1}\frac{1}{m!}\left(\ln\frac{t}{y}\right)^{m}.
\end{align*}

When $0\leq x=y\leq t$, just notice that there is no term $\bigcap_{nx\leq s<ny}
	\left\{\tilde{Z}_{nt}^{l,s+1}=\tilde{Z}_{nt}^{l,s}\right\}$ in the above proof, and the representation for $P_{n,2}$ is the same.

When $0\leq y<x\leq t$, for fixed $1\leq m\leq k_{r}-1$ and sequence $ny=v_{0}\leq v_1<v_2<\cdots<v_m<nt$, define $p:=\max\{0\leq i\leq m:v_{i}<nx\}$, then $\{v_1,\cdots,v_m\}$ can be separated into $\{v_1,\cdots,v_p\}$ and $\{v_{p+1},\cdots,v_m\}$, where $\{v_{1},v_{0}\}=\{v_{m+1},v_{m}\}:=\emptyset$. Then with notation $\prod_{s\in\emptyset}=1$, $P_{n,2}$ can be represented as
\begin{align*}
	&
	\sum_{m=1}^{k_{r}-1}
	\sum_{ny\leq v_1<v_2<\cdots<v_m<nt}
	\left(
	\prod_{\substack{s\in \{v_1,\cdots,v_p\}}}
	\frac{
	\mathbf{P}
	\left(\tilde{Z}_{nt}^{r,s+1}>\tilde{Z}_{nt}^{r,s}\right)
	}
	{
	\mathbf{P}
	\left(\tilde{Z}_{nt}^{r,s+1}=\tilde{Z}_{nt}^{r,s}\right)
	}
	\prod_{s=ny}^{nx-1}
	\mathbf{P}
	\left(\tilde{Z}_{nt}^{r,s+1}=\tilde{Z}_{nt}^{r,s}\right)
	\right.
	\\
	&\left.
	\prod_{s\in \{v_{p+1},\cdots,v_m\}}
	\frac{
	\mathbf{P}
	\left(\tilde{Z}_{nt}^{l,s+1}=\tilde{Z}_{nt}^{l,s},\tilde{Z}_{nt}^{r,s+1}>\tilde{Z}_{nt}^{r,s}\right)
	}
	{
	\mathbf{P}
	\left(\tilde{Z}_{nt}^{l,s+1}=\tilde{Z}_{nt}^{l,s},\tilde{Z}_{nt}^{r,s+1}=\tilde{Z}_{nt}^{r,s}\right)
	}
	\prod_{s=nx}^{nt-1}
	\mathbf{P}
	\left(\tilde{Z}_{nt}^{l,s+1}=\tilde{Z}_{nt}^{l,s},\tilde{Z}_{nt}^{r,s+1}=\tilde{Z}_{nt}^{r,s}\right)
	\right),
\end{align*}
where the second and fourth terms of the product constitute $P_{n,1}$. Thus by \eqref{eq3-3}--\eqref{eq3-5},
	\[\lim_{n\rightarrow\infty}P_{n,2}
	=\lim_{n\rightarrow\infty}P_{n,1} \sum_{m=1}^{k_{r}-1}\sum_{ny\leq v_{1}<\cdots<v_{m}< nt}\left(\frac{1}{v_1v_2\cdots v_m}\right)
	=\frac{xy}{t^{2}(x+1-t)}\sum_{m=1}^{k_{r}-1}\frac{1}{m!}\left(\ln\frac{t}{y}\right)^{m}.\]

So for all $x,y\in[0,t]$,
\begin{equation}\label{eq3-n2}
	\lim_{n\rightarrow\infty}P_{n,2}=\frac{xy}{t^{2}(x+1-t)}\cdot \sum_{m=1}^{k_{r}-1}\frac{1}{m!}\left(\ln\frac{t}{y}\right)^{m}.
\end{equation}

$(iii)$ The limit of $P_{n,3}$:

This part is almost the same as $P_{n,2}$, just notice that by \eqref{eq3-2} and \eqref{eq3-5},
\begin{align*}
	&\lim_{n\rightarrow\infty}\sum_{w=1}^{k_{l}-1}
	\sum_{nx\leq u_1<u_2<\cdots<u_w<nt}
	\prod_{s\in \{u_1,\cdots,u_w\}}
	\frac{
	\mathbf{P}
	\left(\tilde{Z}_{nt}^{l,s+1}>\tilde{Z}_{nt}^{l,s}, \tilde{Z}_{nt}^{r,s+1}=\tilde{Z}_{nt}^{r,s}\right)}
	{\mathbf{P}
	\left(\tilde{Z}_{nt}^{l,s+1}=\tilde{Z}_{nt}^{l,s}, \tilde{Z}_{nt}^{r,s+1}=\tilde{Z}_{nt}^{r,s}\right)}\\
	=&\lim_{n\rightarrow\infty}
	\sum_{w=1}^{k_{l}-1}\sum_{nx\leq u_{1}<u_{2}<\cdots<u_{w}< nt}\left(\frac{1}{u_{1}}-\frac{1}{u_{1}+n-nt}\right)\cdots\left(\frac{1}{u_{w}}-\frac{1}{u_{w}+n-nt}\right)\\
	=&\lim_{n\rightarrow\infty}\sum_{w=1}^{k_{l}-1}\frac{1}{w!}\left(\int_{x}^{t}\left(\frac{1}{s}-\frac{1}{s+1-t}\right)\text{d}s\right)^{w}.
	\end{align*}

Omitting the repetitive proof, we have, for all $x,y\in[0,t]$,
\begin{equation}\label{eq3-n3}
	\lim_{n\rightarrow\infty}P_{n,3}
	=\frac{xy}{t^{2}(x+1-t)}
	\cdot \sum_{w=1}^{k_{l}-1}\frac{1}{w!}\left(\ln\frac{t(x+1-t)}{x}\right)^{w}.
\end{equation}

$(iv)$ The limit of $P_{n,4}$:

Here we only prove the case $0\leq x=y\leq t$, since the other cases are analogues.

For the simplicity of notation, denote
\begin{align*}
	\left\{
	\begin{array}{l}
	A_{s+1}^{1}=\left\{\tilde{Z}_{nt}^{l,s+1}>\tilde{Z}_{nt}^{l,s}, \tilde{Z}_{nt}^{r,s+1}=\tilde{Z}_{nt}^{r,s}\right\},\\
	A_{s+1}^{2}=\left\{\tilde{Z}_{nt}^{l,s+1}=\tilde{Z}_{nt}^{l,s}, \tilde{Z}_{nt}^{r,s+1}>\tilde{Z}_{nt}^{r,s}\right\},\\
	A_{s+1}^{3}=\left\{\tilde{Z}_{nt}^{l,s+1}>\tilde{Z}_{nt}^{l,s}, \tilde{Z}_{nt}^{r,s+1}>\tilde{Z}_{nt}^{r,s}\right\},\\
	A_{s+1}^{4}=\left\{\tilde{Z}_{nt}^{l,s+1}=\tilde{Z}_{nt}^{l,s}, \tilde{Z}_{nt}^{r,s+1}=\tilde{Z}_{nt}^{r,s}\right\}.
	\end{array}
	\right.
\end{align*} 

Then for $1\leq w\leq k_{l}-1$, $1\leq m\leq k_{r}-1$, by \eqref{eq3-1}--\eqref{eq3-5} and converting the summation into integral like $P_{n,2}$ and $P_{n,3}$, we have, as $n\rightarrow\infty$,
\begin{align*}
    &\mathbf{P}
    \left(
    \big|\{nx\leq i\leq nt-1:\tilde{Z}_{nt}^{l,i+1}>\tilde{Z}_{nt}^{l,i}\}\big|=w,
	\big|\{nx\leq i\leq nt-1:\tilde{Z}_{nt}^{r,i+1}>\tilde{Z}_{nt}^{r,i}\}\big|=m
	\right)\\
	=&
	\sum_{\substack{nx\leq u_{1}<\cdots<u_{w}< nt\\nx\leq v_{1}<\cdots<v_{m}< nt}}
	\mathbf{P}
	\left(
	\bigcap_{s\in\{u_{1},\cdots,u_{w}\}}
	\left\{\tilde{Z}_{nt}^{l,s+1}>\tilde{Z}_{nt}^{l,s}\right\}
	\bigcap_{\substack{nx\leq s< nt\\s\notin\{u_{1},\cdots,u_{w}\}}}
	\left\{\tilde{Z}_{nt}^{l,s+1}=\tilde{Z}_{nt}^{l,s}\right\}
	\right.\\
	&\left.\qquad\qquad\qquad\qquad\qquad
	\bigcap_{s\in\{v_{1},\cdots,v_{m}\}}
	\left\{\tilde{Z}_{nt}^{r,s+1}>\tilde{Z}_{nt}^{r,s}\right\}
	\bigcap_{\substack{nx\leq s< nt\\s\notin\{v_{1},\cdots,v_{m}\}}}
	\left\{\tilde{Z}_{nt}^{r,s+1}=\tilde{Z}_{nt}^{r,s}\right\}
	\right)\\
	=&P_{n,1}\cdot
	\sum_{\substack{nx\leq u_{1}<\cdots<u_{w}< nt\\nx\leq v_{1}<\cdots<v_{m}< nt}}
	\prod_{\substack{s\in\{u_{1},\cdots,u_{w}\}\\\backslash\{v_{1},\cdots,v_{m}\}}}
	\frac
	{\mathbf{P}(A_{s+1}^{1})}{\mathbf{P}(A_{s+1}^{4})}
	\prod_{\substack{s\in\{v_{1},\cdots,v_{m}\}\\\backslash\{u_{1},\cdots,u_{w}\}}}
	\frac
	{\mathbf{P}(A_{s+1}^{2})}{\mathbf{P}(A_{s+1}^{4})}
	\prod_{\substack{s\in\{u_{1},\cdots,u_{w}\}\\\cap\{v_{1},\cdots,v_{m}\}}}
	\frac{\mathbf{P}(A_{s+1}^{3})}{\mathbf{P}(A_{s+1}^{4})}
	\\
	=&P_{n,1}\cdot
	\left(
	\sum_{\substack{nx\leq u_{1}<\cdots<u_{w}< nt\\nx\leq v_{1}<\cdots<v_{m}< nt\\\{u_{1},\cdots,u_{w}\}\cap\{v_{1},\cdots,v_{m}\}=\emptyset}}
	\prod_{s\in\{u_{1},\cdots,u_{w}\}}
	\frac
	{\mathbf{P}(A_{s+1}^{1})}{\mathbf{P}(A_{s+1}^{4})}
	\prod_{s\in\{v_{1},\cdots,v_{m}\}}
	\frac
	{\mathbf{P}(A_{s+1}^{2})}{\mathbf{P}(A_{s+1}^{4})}
	\right.
	\\
	+&\left.
	\sum_{\substack{nx\leq u_{1}<\cdots<u_{w}< nt\\nx\leq v_{1}<\cdots<v_{m}< nt\\\{u_{1},\cdots,u_{w}\}\cap\{v_{1},\cdots,v_{m}\}\neq\emptyset}}
	\prod_{\substack{s\in\{u_{1},\cdots,u_{w}\}\\\backslash\{v_{1},\cdots,v_{m}\}}}
	\frac
	{\mathbf{P}(A_{s+1}^{1})}{\mathbf{P}(A_{s+1}^{4})}
	\prod_{\substack{s\in\{v_{1},\cdots,v_{m}\}\\\backslash\{u_{1},\cdots,u_{w}\}}}
	\frac{\mathbf{P}(A_{s+1}^{2})}{\mathbf{P}(A_{s+1}^{4})}
	\prod_{\substack{s\in\{u_{1},\cdots,u_{w}\}\\\cap\{v_{1},\cdots,v_{m}\}}}
	\frac{\mathbf{P}(A_{s+1}^{3})}{\mathbf{P}(A_{s+1}^{4})}
	\right)\\
	\sim& P_{n,1}\cdot \left(\sum_{nx\leq u_{1}<\cdots<u_{w}< nt}
	\prod_{s\in\{u_{1},\cdots,u_{w}\}}
	\frac{\mathbf{P}(A_{s+1}^{1})}{\mathbf{P}(A_{s+1}^{4})}
	\right)\left(\sum_{nx\leq v_{1}<\cdots<v_{m}< nt}
	\prod_{s\in\{v_{1},\cdots,v_{m}\}}
	\frac{\mathbf{P}(A_{s+1}^{2})}{\mathbf{P}(A_{s+1}^{4})}
	\right)\\
	\rightarrow&\frac{xy}{t^{2}(x+1-t)}
	\cdot \frac{1}{w!}\left(\ln\frac{t(x+1-t)}{x}\right)^{w}
	\cdot \frac{1}{m!}\left(\ln\frac{t}{y}\right)^{m}.
	\end{align*}
	
	Thus for all $x,y\in[0,t]$,
\begin{equation}\label{eq3-n4}
	\lim_{n\rightarrow\infty}P_{n,4}
	=\frac{xy}{t^{2}(x+1-t)}
	\left(\sum_{w=1}^{k_{l}-1}\frac{1}{w!}\left(\ln\frac{t(x+1-t)}{x}\right)^{w}\right)
	\left(\sum_{m=1}^{k_{r}-1}\frac{1}{m!}\left(\ln\frac{t}{y}\right)^{m}\right).
\end{equation}
	
And the theorem follows by \eqref{eq3-n1}--\eqref{eq3-n4}.
\end{proof}

Since the generation of most recent common ancestor of all particles at generation $[nt]$ is the minimum of that of the left and right parts at generation $[nt]$, we can easily get,
\begin{proof}[Proof of Corollary \ref{cor:MRCA}]
	Define 
	\[\tilde{G}_{nt}=\min\{0\leq i\leq nt | \tilde{Z}_{nt}^{l,i}=\tilde{Z}_{nt}^{l},\tilde{Z}_{nt}^{r,i}=\tilde{Z}_{nt}^{r}\},\]
	then $\mathscr{L}(G_{nt}|Z_{n}>0)=\mathscr{L}(nt-\tilde{G}_{nt})$. Thus	
	\[
	\lim_{n\rightarrow\infty}\mathbf{P}\left(\frac{G_{nt}}{n}\leq x\Big|Z_{n}>0\right)=\lim_{n\rightarrow\infty}\mathbf{P}\left(\frac{nt-\tilde{G}_{nt}}{n}\leq x\right)=1-\lim_{n\rightarrow\infty}\mathbf{P}\left(\tilde{G}_{nt}\leq n(t-x)\right),\]
	and the corollary follows by Theorem \ref{th:MRCA-in}.
\end{proof}

Using the same method and notations as in Theorem \ref{th:decom-3l} and \ref{th:decom-3r}, we can obtain that for all $n\geq0$ and $0<t<1$, $k\geq1$,
\begin{align}\label{eq:decom-3}
\mathscr{L}\left(\dot{Z}_{nt}^{l},\dot{Z}_{nt}^{r}\right)=\mathscr{L}\left(\frac{a_{N_{nt}^{l,k}T_{nt}^{l,k}}}{a_{nt}}\dot{Z}_{N_{nt}^{l,k}T_{nt}^{l,k}}^{l}+\sum_{i=1}^{k}\frac{b_{N_{nt}^{l,i}T_{nt}^{l,i}}}{a_{nt}}\sum_{j=1}^{\tilde{D}_{nt}^{l,i}}\dot{L}_{N_{nt}^{l,i}T_{nt}^{l,i};j},\right.\nonumber\\
\left.\frac{c_{N_{nt}^{r,k}T_{nt}^{r,k}}}{c_{nt}}\dot{Z}_{N_{nt}^{r,k}T_{nt}^{r,k}}^{r}+\sum_{i=1}^{k}\frac{d_{\tilde{G}_{nt}^{r,i}-1}}{c_{nt}}\sum_{j=1}^{\tilde{D}_{nt}^{r,i}}\dot{R}_{\tilde{G}_{nt}^{r,i}-1;j}\right),
\end{align}
where by the independence of the subtrees, $\{\dot{L}_{n,t;j}, \dot{R}_{nt;j}\}_{n,t,j}$ are independent. Then the asymptotic independence between $\tilde{Z}_{nt}^{l}$ and $\tilde{Z}_{nt}^{r}$ is therefore a direct consequence of the preceding results.
\begin{proof}[Proof of Theorem \ref{th:Z-in}]
		Combining Theorem \ref{th:D-r}, \ref{th:D-l} and taking $k=nt$ in \eqref{eq:decom-3}, we have, as $n\rightarrow\infty$, $(n^{-1}\tilde{Z}_{nt}^{l}$, $n^{-1}\tilde{Z}_{nt}^{r})$ has the same limiting distribution as 
		\[\left(\frac{t(1-t)\sigma^2}{2}\sum_{i=1}^{nt}\frac{b_{N_{nt}^{l,i}T_{nt}^{l,i}}}{a_{nt}}\dot{L}_{N_{nt}^{l,i}T_{nt}^{l,i}},\frac{t\sigma^{2}}{2}\sum_{i=1}^{nt}\frac{d_{\tilde{G}_{nt}^{r,k}-1}}{c_{nt}}\dot{R}_{\tilde{G}_{nt}^{r,i}-1}\right),\]
	where $\{\dot{L}_{n,t}\}$ and $\{\dot{R}_{n}\}$ are independent. Thus by Proposition 8(3) of Davydov and Novikov \cite{Davydov} and  Theorem \ref{th:MRCA-in}, $n^{-1}\tilde{Z}_{nt}^{l}$ and $n^{-1}\tilde{Z}_{nt}^{r}$ are asymptotically independent as $n\rightarrow\infty$.
\end{proof}

\section{Proof of Theorem \ref{th:Spitzer}}\label{section:5}
\begin{proof}[Proof of Theorem \ref{th:Spitzer}]
By  Theorem \ref{th:yaglom-l} and \ref{th:yaglom-r},
	\[n^{-1}\tilde{Z}_{nt}^{l}\overset{\text{d}}{\rightarrow}U_t,\quad n^{-1}\tilde{Z}_{nt}^{r}\overset{\text{d}}{\rightarrow}V_t,\]
	where $U_{t}$ and $V_{t}$ are random variables having exponential distributions with parameters $2/(t(1-t)\sigma^{2})$ and $2/(t\sigma^{2})$ respectively. Thus by Theorem \ref{th:Z-in} and Proposition 1 of Davydov and Novikov \cite{Davydov}, for all $\lambda>0$,
	\[\lim_{n\rightarrow\infty}\mathbf{E}e^{-\lambda(n^{-1}\tilde{Z}_{nt}^{l}+n^{-1}\tilde{Z}_{nt}^{r})}=\lim_{n\rightarrow\infty}\mathbf{E}e^{-\lambda n^{-1}\tilde{Z}_{nt}^{l}}\cdot\lim_{n\rightarrow\infty}\mathbf{E}e^{-\lambda n^{-1}\tilde{Z}_{nt}^{r}}=\mathbf{E}e^{-\lambda U_t}\cdot \mathbf{E}e^{-\lambda V_t}=\mathbf{E}e^{-\lambda(U_t+V_t)},\]
	which implies that $\mathscr{L}\left(n^{-1}Z_{nt}|Z_{n}>0\right)\overset{\text{d}}{\rightarrow}U_t+V_t.$
\end{proof}

\appendix
\section{Appendix}
For the distribution function and mean value of nested uniform random variable sequence, we have the following proposition.
\begin{proposition}\label{prop:nurv}
If $\{U_{k}\}_{k\geq1}$ is a sequence of nested uniform random variables on $[a,b]$, then

(i) for all $k\geq1$, the distribution function of $U_{k}$ is
\[\mathbf{P}\left(U_{k}\leq x\right)=\frac{x-a}{b-a}\sum_{m=0}^{k-1}\frac{1}{m!}\left(\ln\frac{b-a}{x-a}\right)^{m},\quad a\leq x\leq b;\]

(ii) for all $k\geq1$,  the mean value of $U_{k}$ is
\[\mathbf{E}U_{k}=a+\frac{b-a}{2^{k}}.\]
\end{proposition}
\begin{proof}Since for all $a<b$, $\{(b-a)^{-1}(U_{k}-a)\}_{k}$ is a sequence of nested random variables on $[0,1]$, we just need to prove the proposition for $a=0,b=1$.

$(i)$ We will prove it by induction. Obviously, it is true for $k=1$. If it is true for some $k\geq1$, then for $k+1$,
	\begin{align*}
		\mathbf{P}\left(U_{k+1}\leq x\right)
		=&\mathbf{E}\left[\mathbf{P}\left(U_{k+1}\leq x|U_{k}\right)\right]\\
		=&\mathbf{P}\left(U_{k}\leq x\right)+\mathbf{E}\left(\frac{x}{U_{k}};U_{k}> x\right)\\
		=&\mathbf{P}\left(U_{k}\leq x\right)+\int_{x}^{1}\frac{x}{s}\left(\frac{1}{(k-1)!}\left(\ln\frac{1}{s}\right)^{k-1}\right)ds\\
		=&\mathbf{P}\left(U_{k}\leq x\right)+\frac{x}{k!}\left(\ln\frac{1}{x}\right)^{k}\\
		=&x\sum_{m=0}^{k}\frac{1}{m!}\left(\ln\frac{1}{x}\right)^{m}.
	\end{align*}
	
$(ii)$ Just notice that $\mathbf{E}U_{1}=2^{-1}$ and $\mathbf{E}U_{k}=\mathbf{E}[\mathbf{E}(U_{k}|U_{k-1})]=2^{-1}\mathbf{E}U_{k-1}$.
\end{proof}

\begin{lemma}\label{lem:Lnt}
Let $L_{nt}$ be the random variable with distribution $\mathscr{L}(L_{nt})=\mathscr{L}(Z_{nt}|Z_{nt}>0,Z_{n}=0)$. Suppose $\alpha=1$ and $\sigma^{2}<\infty$, then for any fixed $0<t<1$,
\begin{align*}
	&\text{(i)}\lim_{n\rightarrow\infty}\mathbf{P}\left(\frac{L_{nt}}{n}\geq x\right)=\exp\left(-\frac{2x}{t(1-t)\sigma^{2}}\right), \quad x\geq 0.\\
	&\text{(ii)}\lim_{n\rightarrow\infty}n^{-1}\mathbf{E}L_{nt}=\frac{t(1-t)\sigma^2}{2}.\\
	&\text{(iii)}\lim_{n\rightarrow\infty}n^{-2}\mathbf{E}L_{nt}^{2}=\frac{[t(1-t)\sigma^2]^2}{2}.
\end{align*}	
\end{lemma}
\begin{proof}
Observe that for $j\geq0$,
\[\mathbf{P}(Z_{nt}=j|Z_{nt}>0,Z_{n}=0)
		=\frac{\mathbf{P}(Z_{nt}>0)\mathbf{P}(Z_{nt}=j|Z_{nt}>0)[\mathbf{P}(Z_{n(1-t)}=0)]^{j}}{\mathbf{P}(Z_{n}=0)-\mathbf{P}(Z_{nt}=0)}.\]
(i) For  $\forall a>0$,
	\begin{align*}
		\mathbf{E}\left(e^{-a\frac{Z_{nt}}{n}}\big|Z_{nt}>0, Z_{n}=0\right)
		=&\sum_{j=0}^{\infty}e^{-\frac{aj}{n}}\mathbf{P}(Z_{nt}=j|Z_{nt}>0,Z_{n}=0)\\
		=&\frac{\mathbf{P}(Z_{nt}>0)}{\mathbf{P}(Z_{n}=0)-\mathbf{P}(Z_{nt}=0)}\mathbf{E}\left(e^{-a\frac{Z_{nt}}{n}}[f_{n(1-t)}(0)]^{Z_{nt}}\big|Z_{nt}>0\right)\\
		=&\left(\frac{1}{1-t}+o(1)\right)\mathbf{E}\left(e^{-\frac{Z_{nt}}{nt}(at-nt\log f_{n(1-t)}(0))}\big|Z_{nt}>0\right).
	\end{align*}
	
	But
	$\lim_{n\rightarrow\infty}(at-nt\log f_{n(1-t)}(0))=at+\left(\frac{t}{1-t}\right)\frac{2}{\sigma^{2}},$
	and by Yaglom's Theorem,
	\[\lim_{n\rightarrow\infty}\mathbf{E}\left(e^{-a\frac{Z_{nt}}{nt}}\big|Z_{nt}>0\right)=\frac{1}{1+a\frac{\sigma^{2}}{2}},\]
	thus
	\[\lim_{n\rightarrow\infty}\mathbf{E}\left(e^{-a\frac{Z_{nt}}{n}}\big|Z_{nt}>0, Z_{n}=0\right)=\frac{1}{1-t}\cdot \frac{1}{1+a t\frac{\sigma^{2}}{2}+\frac{t}{1-t}}=\left(1+at(1-t)\frac{\sigma^2}{2}\right)^{-1},\]
	and the conclusion follows.
	
(ii) Just notice that 
\begin{align*}
	\mathbf{E}\left(Z_{nt}\big|Z_{nt}>0, Z_{n}=0\right)
	&=\frac{\sum_{j=1}^{\infty}j\mathbf{P}(Z_{nt}=j)[\mathbf{P}(Z_{n(1-t)}=0)]^{j}}{\mathbf{P}(Z_{n}=0)-\mathbf{P}(Z_{nt}=0)}\\
	&=\frac{\mathbf{P}(Z_{n}=0)\mathbf{E}\left(Z_{nt}\big|Z_{n}=0\right)}{\mathbf{P}(Z_{n}=0)-\mathbf{P}(Z_{nt}=0)},
\end{align*}
	and by Proposition \ref{prop:moment-l1} with $i=nt$, $\mathbf{E}\left(Z_{nt}\big|Z_{n}=0\right)\rightarrow(1-t)^2.$

(iii) Similarly,
\[\mathbf{E}\left(Z_{nt}^2\big|Z_{nt}>0, Z_{n}=0\right)
	=\frac{\mathbf{P}(Z_{n}=0)\mathbf{E}\left(Z_{nt}^2\big|Z_{n}=0\right)}{\mathbf{P}(Z_{n}=0)-\mathbf{P}(Z_{nt}=0)},\]
and by Proposition \ref{prop:moment-l2} with $i=nt$, $\mathbf{E}\left(Z_{nt}^2\big|Z_{n}=0\right)\sim(1-t)^3\sigma^{2}nt.$
\end{proof}

\end{document}